\documentclass{article}

\usepackage{authblk}
\usepackage[utf8]{inputenc}
\usepackage[T1]{fontenc}
\usepackage{pslatex}
\usepackage[english]{babel}
\usepackage{fancyhdr}
\usepackage{hyperref}
\usepackage{graphicx}
\usepackage{amsmath}
\usepackage{amsthm}
\usepackage{amsfonts}
\usepackage{mathtools}
\usepackage{amssymb}
\usepackage{algorithm}
\usepackage{algorithmic}
\usepackage{microtype}
\usepackage{rotating}
\usepackage{todonotes}
\usepackage{tikz}
\usepackage{tikz-qtree}
\usepackage{xparse}
\usepackage{appendix}
\usepackage{url}
\usepackage{tabularx}
\usepackage{paralist}
\usepackage{booktabs}
\usepackage{multirow}
\usepackage{subcaption}
\usepackage{cancel}
\usetikzlibrary{calc,shapes.callouts,shapes.arrows}
\usetikzlibrary{automata}
\usetikzlibrary{positioning}
\usetikzlibrary{decorations.pathreplacing}
\usetikzlibrary{decorations.pathmorphing}
\usetikzlibrary{arrows}
\usetikzlibrary{shapes}

\tikzset{
	mybrace/.style={decorate,decoration={brace,aspect=#1}}
}

\newcommand{\N}{\mathbb{N}}
\newcommand{\Z}{\mathbb{Z}}
\newcommand{\F}{\mathbb{F}}

\newtheorem{definition}{Definition}
\newtheorem{theorem}{Theorem}
\newtheorem{lemma}{Lemma}

\providecommand{\keywords}[1]{\textbf{\textit{Keywords }} #1}

\begin{document}

\title{Combinatorial Designs and Cellular Automata: A Survey}

\author[1]{Luca Manzoni}
\author[2]{Luca Mariot}
\author[3]{Giuliamaria Menara}

\affil[1 ]{{\normalsize \normalsize Department of Mathematics, Informatics and Geosciences, University of Trieste, Via Valerio 12/1, 34127 Trieste, Italy} \\
	
	{\small \texttt{lmanzoni@units.it}}}

\affil[2 ]{{\normalsize Semantics, Cybersecurity and Services Group, University of Twente, 7522 NB Enschede, The Netherlands} \\
	
	{\small \texttt{l.mariot@utwente.nl}}}

\affil[3 ]{{\normalsize Department of Informatics, Systems and Communication, University of Milano-Bicocca, Viale Sarca 336, 20126 Milano, Italy} \\
	
	{\small \texttt{giuliamaria.menara@unimib.it}}}
	
\maketitle

\begin{abstract}
Cellular Automata (CA) are commonly investigated as a particular type of dynamical systems, defined by shift-invariant local rules. In this paper, we consider instead CA as algebraic systems, focusing on the combinatorial designs induced by their short-term behavior. Specifically, we review the main results published in the literature concerning the construction of mutually orthogonal Latin squares via bipermutive CA, considering both the linear and nonlinear cases. We then survey some significant applications of these results to cryptography, and conclude with a discussion of open problems to be addressed in future research on CA-based combinatorial designs.
\end{abstract}


\section{Introduction}
\label{sec:intro}
Among all natural computing models, Cellular Automata (CA) stand out for their simplicity and massive parallelism. Indeed, a CA can be broadly defined as a \emph{shift-invariant} transformation over a regular lattice of \emph{cells}, whose states are updated according to the synchronous application of a \emph{local} rule. The two properties of shift-invariance and locality are powerful enough to enable the modelization of a wide variety of dynamical systems. Indeed, since their inception in the 50s by Ulam~\cite{ulam52} and Von Neumann~\cite{vonneumann66}, CA have  been used to simulate a plethora of phenomena in physics~\cite{toffoli84,lent94,wolf04}, biology~\cite{ermentrout93,sirakoulis03,deutsch21}, economics~\cite{goldenberg01,leydesdorff02,caruso07}, ecology~\cite{balzter98,rohde05,garciaduro18} and other fields.

Beyond simulating complex systems, researchers successfully employed CA to implement efficient algorithms and computational devices. The case of \emph{cryptography} is perhaps one of the most prominent: since  Wolfram's work on CA-based pseudorandom number generators (PRNG)~\cite{wolfram85}, the literature flourished with many works leveraging CA to design the most disparate cryptographic primitives. Aside from PRNG and stream ciphers, examples include block ciphers~\cite{daemen94,szaban08,picek17}, hash functions~\cite{daemen91,mihaljevic98,keccak}, public-key encryption schemes~\cite{kari92,clarridge09,applebaum10} and secret sharing schemes~\cite{rey05,ml-acri-2014,ml-acri-2018}. The main motivations to use CA in cryptography are their inherent parallelism (which lends itself to very efficient hardware implementations) and the possibility of having quite complex and unpredictable dynamics starting from simple rules. CA can thus be used to realize transformations with good \emph{confusion} and \emph{diffusion}, two fundamental properties set forth by Shannon~\cite{shannon49} that every secure cryptosystem should satisfy to withstand cryptanalytic attacks.

However, it turned out that considering only the dynamic behavior of certain CA as a surrogate for confusion and diffusion is not a sound approach. For instance, Wolfram's CA-based PRNG has been proven vulnerable to a few attacks~\cite{meier-staff,koc-ca}, despite the fact that the underlying local rule induces a chaotic CA. The reason lies in the observation that most of the properties relevant for cryptography have an algebraic or combinatorial nature, rather than a system-theoretic one. In particular, focusing on the algebraic properties that arise from the short-term behavior of CA proved to be a successful approach for the design of several cryptographic primitives~\cite{daemen94,keccak,mariot19}.

Among the various works published in this area, a recent thread considers the perspective of combinatorial designs induced by CA when interpreted as algebraic systems. Although the first works that investigated this approach date back to the 90s~\cite{p-complexsyst-1992,md-complexsyst-1996,m-physicad-1998}, the explicit connection to cryptographic applications has been considered only in the last few years~\cite{mgfl-desi-2020,gmp-desi-2023}.

The aim of this survey is to systematically review the main results in the literature concerning the construction of combinatorial designs through CA, particularly emphasizing the case of Mutually Orthogonal Latin Squares (MOLS). Accordingly, we review the main theoretical results concerning the construction of MOLS via bipermutive CA both in the linear and nonlinear case, a research thread that has been mainly initiated and developed by the second author of this manuscript. Further, this survey gives an overview of the applications of CA-based combinatorial designs to cryptography and coding theory, including the construction of:

\begin{itemize}
    \item Threshold secret sharing schemes~\cite{mgfl-desi-2020,ml-acri-2018}.
    \item Pseudorandom number generators~\cite{m-naco-2023}.
    \item Hadamard matrices and bent Boolean functions~\cite{gmp-iacr-2020,gmp-desi-2023}.
    \item Binary orthogonal arrays and correlation-immune Boolean functions~\cite{mariot2023building}.
\end{itemize}

In surveying these results and applications, this work also identifies some interesting open problems for future research about the construction of combinatorial designs based on CA.

The rest of this paper is organized as follows. Section~\ref{sec:bg} covers the background concepts related to cellular automata and combinatorial designs used throughout this work. Section~\ref{sec:ca-alg} introduces the interpretation of CA as algebraic systems, and gives a survey of the first works that adopted it to investigate questions not directly related to cryptography. Then, the section discusses the basic motivation for CA-based combinatorial designs in the design of cryptographic applications, showing the basic construction of Latin squares from bipermutive CA. Next, Section~\ref{sec:mols-ca} focuses on \emph{Mutually Orthogonal Cellular Automata} (MOCA), which are families of bipermutive CA that form \emph{Mutually Orthogonal Latin Squares} (MOLS). The section covers both the well-developed theory of linear MOCA, and the few results discovered so far for the nonlinear case, which prompt several open questions. Section~\ref{sec:appl} surveys some applications of MOCA in cryptography. Finally, Section~\ref{sec:outro} discusses a few interesting open problems and directions for future research concerning the CA perspective on combinatorial designs.

\section{Background}
\label{sec:bg}
In this section, we give the preliminary definitions and results needed for the rest of the paper. We start in Section~\ref{subsec:ca} by introducing a CA model that can be interpreted as an algebraic system. In Section~\ref{subsec:cd} we define the main combinatorial designs of our interest, namely mutually orthogonal Latin squares and orthogonal arrays.

\subsection{Cellular Automata}
\label{subsec:ca}

Cellular Automata (CA) represent a specific category of parallel computational models. Informally, a cellular automaton consists of a grid of \emph{cells}, with each cell employing a \emph{local rule} to determine its next state based on its \emph{neighborhood}. The overall state of the CA, known as the \emph{global state}, encompasses the configuration of all the cells within the grid at a specific moment. The dynamic behavior of the CA unfolds as all cells update simultaneously during discrete time steps, something which can be formalized by defining a \emph{global rule} in terms of the local rule.

There are various ways to study the properties of CA. One of the most common approaches is to consider a CA on an \emph{infinite} one-dimensional grid, with the state of each cell ranging over a finite alphabet $\Sigma$. The set of all possible global states, denoted as $\Sigma^\Z$, can then be endowed with the structure of a compact metric space under the \emph{Cantor distance}. This construction results in the so-called \emph{full-shift space} studied in the field of symbolic dynamics~\cite{lind21}. The \emph{Curtis-Hedlund-Lyndon} theorem~\cite{hedlund69} characterizes CA as those endomorphisms $F: \Sigma^\Z \to \Sigma^\Z$ of the full-shift space that are both uniformly continuous with respect to the Cantor distance and \emph{shift-invariant}. This means that a CA commutes with the \emph{shift operator} $\sigma: \Sigma^\Z \to \Sigma^\Z$ that moves by one place to the left each element of a bi-infinite string.

The topological view is quite powerful to characterize the long-term dynamics of infinite CA. However, in real-world applications---such as cryptography---one is naturally constrained to use finite arrays only, since the CA has to be implemented on a physical device or programmed in software with limited memory. This introduces the problem of updating the boundary cells of the CA, as they lack sufficient neighbors for the evaluation of the local rule. Various workarounds exist in the literature to address this issue. An often-used solution is to consider a finite one-dimensional array of $n$ cells with \emph{periodic boundary conditions}, where the last (rightmost) cell precedes the first (leftmost) cell and vice versa. In this way, each cell has a complete neighborhood and can thus evaluate the local rule to compute its next state. The whole CA can then be iterated indefinitely, although the long-term dynamics of the resulting system is of course ultimately periodic: if the CA is composed of $n$ cells, the dynamic evolution will enter a cycle after at most $|\Sigma|^n$ evaluations of the global rule. Thus, a periodic CA is iterated for a number of steps that is polynomial in the size of the array. To a certain extent, the topological perspective can still be used in this context, since periodic CA correspond to the subset of \emph{spatially periodic configurations} of the full-shift space~\cite{kari12}.

Periodic CA have been thoroughly investigated as finite dynamical systems to implement cryptographic primitives such as pseudorandom number generators~\cite{wolfram85,formenti14,leporati13}, block ciphers~\cite{gutowitz93,marconi06,szaban08} and S-boxes~\cite{daemen94,keccak,mariot19}. We do not delve further into the details of this research thread, as periodic CA are not the focus of this survey. The interested reader may find additional information on related works in~\cite{mariot22,mariot24}.

As we mentioned in the Introduction, the topological dynamics viewpoint is often misleading for cryptographic applications that leverage CA. As a matter of fact, many security criteria are often formalized in terms of algebraic properties, rather than as system-theoretic ones. Furthermore, in several cryptographic applications one does not even need to iterate functions for multiple time steps: often, a single application of a mapping is enough, since it is usually combined with other types of operations. This motivates the following CA model, which we consider in the rest of this paper:

\begin{definition}
    \label{def:nbca}
    Let $\Sigma$ be a finite alphabet and $n,d \in \N$ with $n \geq d$. 
    Additionally, let $f \in \Sigma^d \to \Sigma$ be a local rule of $d$ variables.
    The \emph{No Boundary Cellular Automaton (NBCA)} $F: \Sigma^n \to \Sigma^{n-d+1}$ is the function defined for all $x\in \Sigma^n$ as:
    \[
    F(x_1,\dots,x_n) = (f(x_1,\dots,x_d),f(x_2,\dots,x_{d+1}),\dots,f(x_{n-d+1},\dots,x_n)).
    \]
\end{definition}
In other words, a NBCA (in the following, simply a CA) is defined as a vectorial function where each coordinate function $F_i: \Sigma^n \to \Sigma$ evaluates the local rule $f$ on the neighborhood formed by the $i$-th input cell and the $i+d-1$ cells to its right. Since there are $n-d+1$ output cells in total, the neighborhood is always complete. The parameter $d$ is also called the \emph{diameter} of the CA.

Several choices can be made for the alphabet, with the easiest being the \emph{binary} alphabet $\Sigma = \{0,1\}$. In what follows, we also consider the more general case where $\Sigma = \F_q$ is the \emph{finite field} of order $q$, where $q$ is a power of a prime number~\cite{lidl94}. Clearly, when $q=2$ we have $\F_2 = \{0,1\}$ which coincides again with the binary case. The sum and multiplication operations on $\F_2$ correspond respectively to the XOR (denoted as $\oplus$) and to the logical AND (denoted by concatenation of the operands) of two bits $a,b \in \F_2$.

Once the alphabet is fixed, there are different ways to represent the local rule of a CA. The most straightforward one is the \emph{rule table}, which enumerates the next state of a cell for each possible neighborhood combination. In the binary case, the local rule $f: \F_2^d \to \F_2$ is basically a \emph{Boolean function} of $d$ variables; thus, the rule table is the \emph{truth table} representation of $f$. If a total order is fixed on the input vectors of $\F_2^d$ (e.g.,  the lexicographic order), the truth table of $f$ can be simply identified by its output column $f(x)$, which is a $2^d$-bit string. In the CA terminology, the decimal encoding of $f(x)$ is also known as the \emph{Wolfram code} of the local rule~\cite{wolfram83}.

\begin{figure}[t]
    \centering
    \begin{minipage}{0.3\textwidth}
        \centering
        \resizebox{3cm}{!}{%
            \centering
            \begin{tikzpicture}
                [->,auto,node distance=1.5cm, empt node/.style={font=\sffamily,inner
                    sep=0pt}, rect
                node/.style={rectangle,draw,font=\bfseries,minimum size=0.5cm, inner
                    sep=0pt, outer sep=0pt}]
                
                \node [empt node] (c)   {};
                \node [rect node] (c1) [right=0.1cm of c] {$1$};
                \node [rect node] (c2) [right=0cm of c1] {$0$};
                \node [rect node] (c3) [right=0cm of c2] {$0$};
                \node [rect node] (c4) [right=0cm of c3] {$1$};
                
                \node [empt node] (f1) [above=0.2cm of c3] {{\footnotesize $f(1,0,0) = 1$}};
                
                
                \node [rect node] (p2) [above=0.85cm of c1] {$0$};
                \node [rect node] (p1) [left=0cm of p2] {$1$};
                \node [rect node] (p3) [right=0cm of p2] {$0$};
                \node [rect node] (p4) [right=0cm of p3] {$0$};
                \node [rect node] (p5) [right=0cm of p4] {$0$};
                \node [rect node] (p6) [right=0cm of p5] {$1$};
                
                \node [empt node] (p7) [below=0.2cm of p1] {};
                \node [empt node] (p8) [right=0.07cm of p7] {};
                \node [empt node] (p12) [above=0.5cm of p1.east] {};
                \node [empt node] (p13) [above=0.5cm of p5.east] {};
                \node [empt node] (p14) [above=0.3cm of p13] {\phantom{M}};
                
                \draw [-, mybrace=0.25, decorate, decoration={brace,mirror,amplitude=5pt,raise=0.3cm}]
                (p1.west) -- (p3.east) node [midway,yshift=-0.3cm] {};
                \draw[->] (p8) -- (c1.north);
            \end{tikzpicture}
        }
    \end{minipage}%
    \begin{minipage}{0.4\textwidth}
        \centering
        \begin{tabular}{ccc|c}
            \hline
            $x_i$ & $x_{i+1}$ & $x_{i+2}$ & $f(x)$ \\
            \hline
            0 & 0 & 0 & 0 \\
            0 & 0 & 1 & 1 \\
            0 & 1 & 0 & 1 \\
            0 & 1 & 1 & 0 \\
            1 & 0 & 0 & 1 \\
            1 & 0 & 1 & 0 \\
            1 & 1 & 0 & 0 \\
            1 & 1 & 1 & 1 \\
            \hline
        \end{tabular}
    \end{minipage}%
    \begin{minipage}{0.3\textwidth}
        \centering
        \resizebox{!}{5cm}{%
            \begin{tikzpicture}
                [->,auto,node distance=1.5cm, every loop/.style={min distance=12mm},
                empt node/.style={font=\sffamily,inner sep=0pt,outer sep=0pt},
                circ node/.style={circle,thick,draw,font=\sffamily\bfseries,minimum
                    width=0.8cm, inner sep=0pt, outer sep=0pt}]
                
                \node [empt node] (e1) {};
                \node [circ node] (n00) [above=1.75cm of e1] {$00$};
                \node [circ node] (n01) [right=1.75cm of e1] {$01$};
                \node [circ node] (n10) [left=1.75cm of e1] {$10$};
                \node [circ node] (n11) [below=1.75cm of e1] {$11$};       
                
                \draw [->, thick, shorten >=0pt,shorten <=0pt,>=stealth] (n00) 
                edge[bend left=20] node (f5) [above right]{$1$} (n01);
                \draw [->, thick, shorten >=0pt,shorten <=0pt,>=stealth] (n01)
                edge[bend left=20] node (f5) [below right]{$0$} (n11);
                \draw [->, thick, shorten >=0pt,shorten <=0pt,>=stealth] (n11)
                edge[bend left=20] node (f5) [below left]{$0$} (n10);
                \draw [->, thick, shorten >=0pt,shorten <=0pt,>=stealth] (n10)
                edge[bend left=20] node (f5) [above left]{$1$} (n00);
                \draw[->, thick, shorten >=0pt,shorten <=0pt,>=stealth] (n10) edge[bend
                left=20] node (f1) [above]{$0$} (n01);
                \draw[->, thick, shorten >=0pt,shorten <=0pt,>=stealth] (n01)
                edge[bend left=20] node (f2) [below]{$1$} (n10);
                \draw[->, thick, shorten >=0pt,shorten <=0pt,>=stealth] (n00) edge[loop
                above] node (f3) [above]{$0$} ();
                \draw[->, thick, shorten >=0pt,shorten <=0pt,>=stealth] (n11) edge[loop
                below] node (f4) [below]{$1$} ();
            \end{tikzpicture}
        }
    \end{minipage}
    \caption{Example of NBCA defined by rule 150, together with its truth table and its de Bruijn graph representations.}
    \label{fig:db-150}
\end{figure}
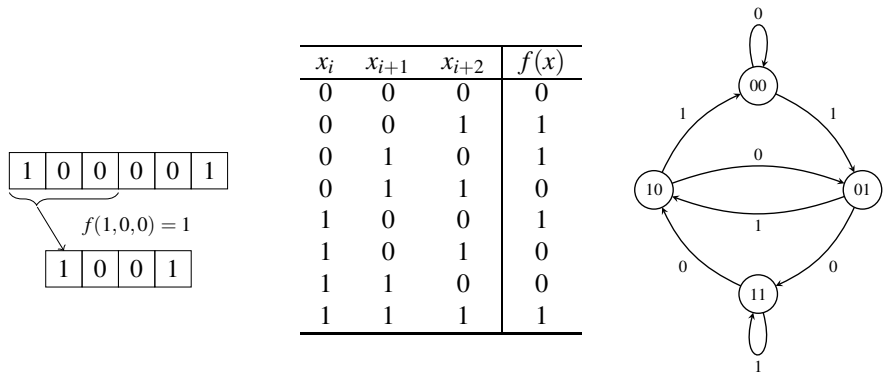

Another common representation for the local rule of a CA is the \emph{de Bruijn graph} $G = (V, E)$ where the set of vertices is $V=\Sigma^{d-1}$, i.e. the set of all possible blocks of $d-1$ cells. Two vertices $u,v \in \Sigma^{d-1}$ are connected by an edge if and only if they \emph{overlap} respectively on the rightmost and leftmost $d-1$ cells, that is, when $u = u_1t$ and $v=tv_1$, with $u_1, v_1 \in \Sigma$ and $t \in \Sigma^{d-2}$. In this case, we can also write the string $x \in \Sigma^d$ of length $d$ as the \emph{fusion} of $u$ and $v$, denoted as $x = u \odot v = u_1tv_1$~\cite{sutner91}. A local rule $f: \Sigma^d \to \Sigma$ can then be represented as a \emph{labeling function} $l: E \to \Sigma$ on the edges of the de Bruijn graph, where $l(u,v) = f(u \odot v)$ for all $(u,v) \in E$. Figure~\ref{fig:db-150} depicts an example of binary CA $F: \F_2^6 \rightarrow \F_2^4$ induced by the local rule $f(x_i, x_{i+1}, x_{i+2}) = x_i \oplus x_{i+1} \oplus x_{i+2}$, along with the truth table and the de Bruijn graph representation of the rule. The Wolfram code of this rule is 150, since it is the decimal encoding of the output column $10010110$ (read in bottom-up order).

An interesting consequence of the de Bruijn graph representation is that the CA input vector corresponds to a \emph{path on the vertices}, merged together with the fusion operator. The output vector of the CA is instead the corresponding \emph{path on the edges} of the graph. For example, the input configuration in Figure~\ref{fig:db-150} can be read as the sequence of vertices $10 \odot 00 \odot 00 \odot 00 \odot 01 = 100001$. The corresponding labels on the edges then give the output configuration $1001$.

A property of local rules that we consider in this paper is permutivity, due to its connection with combinatorial designs. Formally, a local rule $f: \Sigma^d \to \Sigma$ is called \emph{permutive} in the $i$-th variable (for $i \in \{1,\ldots,d\}$) if, by fixing all input coordinates to any value except the $i$-th one, the resulting restriction of $f$ as a function of one variable is a permutation over $\Sigma$. A local rule $f$ that is permutive in the first (respectively, $d$-th) variable is also called \emph{leftmost} (respectively, \emph{rightmost}) permutive. If $f$ is both leftmost and rightmost permutive, then it is called a \emph{bipermutive local rule}. Rule 150 is an example of bipermutive local rule. In fact, any bipermutive local rule over the binary alphabet is a Boolean function $f: \F_2^d \to \F_2$ of the form:
\begin{equation}
    \label{eq:gen-func}
    f(x_1, \ldots, x_d) = x_1 \oplus g(x_2,\ldots,x_{d-1}) \oplus x_d
\end{equation}
for all $x = (x_1,\ldots,x_d) \in \F_2^d$, where $g: \F_2^{d-2} \to \F_2$ is any Boolean function of $d-2$ variables. Hence, the output of $f$ is determined by computing the XOR of the leftmost and rightmost variables, together with a function of the $d-2$ central variables.

CA defined by permutive rules have been thoroughly investigated in the topological dynamics literature. Gilman established in \cite{gilman2006periodic} that infinite one-dimensional cellular automata exhibiting a (bi)permutive behavior are topologically conjugate to one-sided shifts, provided they have an appropriate finite state space. Consequently, it follows that infinite CA demonstrating (bi)permutive characteristics exhibit chaotic behavior as dynamical systems. It is worth noting that this result has been independently rediscovered on multiple occasions, e.g. \cite{cattaneo2000investigating, fagnani1998expansivity, kleveland1997mixing, shereshevsky1992bipermutative}. The connection between the chaotic dynamics of permutive CA and the cryptographic properties related to pseudorandom number generators has been investigated in~\cite{martin08,formenti14,leporati13,leporati14}.

When the alphabet is the finite field $\F_q$, it is possible to introduce the notion of linearity in CA. Specifically, a local rule $f: \F_q^d \to \F_q$ is called \emph{linear} if it is defined as a linear combination of the neighborhood cells, or formally if there exist $a_1,\ldots, a_d \in \F_q$ such that:
\begin{displaymath}
    f(x_1,\ldots,x_d) = a_1x_1 + \ldots + a_dx_d
\end{displaymath}
for all $x = (x_1,\ldots,x_d) \in \F_q^d$, where sum and multiplication correspond to the field operations of $\F_q$. A CA $F: \F_q^n \to \F_q^{n-d+1}$ equipped with a linear local rule $f: \F_q^d \to \F_q$ is a linear map from the vector space $\F_q^n$ to $\F_q^{n-d+1}$, and the global rule is defined as the matrix-vector multiplication $F(x) = M_F \cdot x^\top$ for all $x \in \F_q^n$, where $M_F$ is a \emph{transition matrix} with the following structure:
\begin{equation}
    \label{eq:ca-matr-f}
    M_F = 
    \begin{pmatrix}
        a_1    & \ldots & a_{d-1} & a_d & 0 & \ldots & \ldots & \ldots & \ldots & 0 \\
        0      & a_1    & \ldots  & a_{d-1} & a_d & 0 & \ldots & \ldots & \ldots & 0 \\
        \vdots & \vdots & \vdots & \ddots  & \vdots & \vdots & \vdots & \ddots & \vdots & \vdots \\
        0 & \ldots & \ldots & \ldots & \ldots & 0 & a_1 & \ldots & a_{d-1} & a_d \\
    \end{pmatrix} \enspace .
\end{equation}
As remarked in \cite{mariot2018cryptographic}, the matrix $M_F$ in Equation \eqref{eq:ca-matr-f} is the generator matrix of a cyclic code.
Hence, one can naturally define the polynomial $p_f(X) \in \mathbb{F}_q[X]$ associated to a linear CA F as the generator polynomial of degree at most $d-1$ of the corresponding cyclic code:

\begin{displaymath}
    p_f(X) = a_1 + a_2X + a_3 X^2 + \ldots + a_dX^{d-1} \enspace .
\end{displaymath}
Stated otherwise, the associated polynomial is defined by using the coefficients $a_1, \ldots, a_d$ of the local rule as the coefficients of the increasing powers of the unknown $X$. Notice that a linear CA is bipermutive if and only if both $a_1$ and $a_d$ are not null. In this case, the associated polynomial is of degree $d-1$ and has a nonzero constant term.

\subsection{Combinatorial Designs}
\label{subsec:cd}
Broadly speaking, the main object of interest in the theory of \emph{combinatorial designs} are families of subsets of a finite set that satisfy certain balancedness constraints. Research in this field took hold in the 20th century, although earlier results can be traced back to Euler~\cite{euler-cd}. Nowadays combinatorial designs found applications in the most diverse domains of science, ranging from the design of experiments to cryptography and error-correcting codes.

The scope of combinatorial designs is too large to be covered in this paper. For this reason, here we focus only on the two main breeds of combinatorial designs that have been studied in the context of CA, namely Latin squares and orthogonal arrays. The interested reader can find a more comprehensive treatment of combinatorial designs in standard references such as~\cite{stinson-cd,handbook-cd,denes-ls,hedayat12}.

In the following, we denote by $[N] = \{1,\ldots, N\}$ the set of the first $N$ positive integer numbers, for all $N \in \N$. A Latin square is formally defined as follows:
\begin{definition}
    \label{def:ls}
    A \emph{Latin square} of order $N \in \N$ is a $N \times N$ square matrix $L$ with entries from $[N]$, such that $L(i,j) \neq	L(i,k)$ and $L(j,i) \neq L(k,i)$ for all $i, j, k \in [N]$ .
\end{definition}
Intuitively, an $N\times N$ matrix is a Latin square of order $N$ if and only if each number from $1$ to $N$ occurs exactly once in each row and each column, or equivalently each row and each column is a permutation of $[N]$.

A Latin square of any order $N \in \N$ can be easily constructed by taking the elements of $[N]$ in increasing order for the first row, and then defining the subsequent rows by repeatedly applying cyclic shifts.

The algebraic structure related to the notion of Latin square is the \emph{quasigroup}, for which we give the following definition:
\begin{definition}
    A \emph{quasigroup} of order $N \in \N$ is a pair $(X, \circ)$ where $X$ is a finite set of cardinality $N$, and $\circ$ is a binary operation over $X$ such that for all $x, y \in X$ the two equations $x \circ z = y$ and $z	\circ x = y$ admit a unique solution for all $z \in X$.
\end{definition}
In particular, a finite algebraic structure $(X, \circ)$ is a quasigroup if and only if its \emph{Cayley table} is a Latin square~\cite{stinson-cd}.

We now introduce the property of orthogonality, which takes into account two Latin squares of the same order:
\begin{definition}
    \label{def:ols}
    Two Latin squares $L_1, L_2$ of order $N$ are \emph{orthogonal} if
    \begin{equation}
        (L_1(i_1,j_1),L_2(i_1,j_1)) \neq (L_1(i_2,j_2),L_2(i_2,j_2))
    \end{equation}
    for all distinct pairs of coordinates $(i_1,j_1), (i_2,j_2) \in [N]\times [N]$. 
\end{definition}
Equivalently to the above definition, $L_1$ and $L_2$ are orthogonal if, by \emph{superposing} one on top of the other, all ordered pairs of the Cartesian product $[N] \times [N]$ occur exactly once. A set of $k$ \emph{Mutually Orthogonal Latin Squares} ($k$-MOLS) is a set of Latin squares of the same order that are pairwise orthogonal. As an example, Figure~\ref{fig:ols-ex} depicts a pair of orthogonal Latin squares of order $4$.
\begin{figure}
    \centering
    \begin{minipage}{.3\textwidth}
        \begin{tikzpicture}
            [->,auto,node distance=1.5cm,
            empt node/.style={font=\sffamily,inner sep=0pt,minimum size=0pt},
            rect node/.style={rectangle,draw,font=\sffamily,minimum size=0.8cm, inner sep=0pt, outer sep=0pt}]
            
            \node [rect node] (s11) {$1$};
            \node [rect node] (s12) [right=0cm of s11] {$3$};
            \node [rect node] (s13) [right=0cm of s12] {$4$};
            \node [rect node] (s14) [right=0cm of s13] {$2$};
            
            \node [rect node] (s21) [below=0cm of s11] {$4$};
            \node [rect node] (s22) [right=0cm of s21] {$2$};
            \node [rect node] (s23) [right=0cm of s22] {$1$};
            \node [rect node] (s24) [right=0cm of s23] {$3$};
            
            \node [rect node] (s31) [below=0cm of s21] {$2$};
            \node [rect node] (s32) [right=0cm of s31] {$4$};
            \node [rect node] (s33) [right=0cm of s32] {$3$};
            \node [rect node] (s34) [right=0cm of s33] {$1$};
            
            \node [rect node] (s41) [below=0cm of s31] {$3$};
            \node [rect node] (s42) [right=0cm of s41] {$1$};
            \node [rect node] (s43) [right=0cm of s42] {$2$};
            \node [rect node] (s44) [right=0cm of s43] {$4$};
            
        \end{tikzpicture}
    \end{minipage}%
    \begin{minipage}{.3\textwidth}
        \begin{tikzpicture}
            [->,auto,node distance=1.5cm,
            empt node/.style={font=\sffamily,inner sep=0pt,minimum size=0pt},
            rect node/.style={rectangle,draw,font=\sffamily,minimum size=0.8cm, inner sep=0pt, outer sep=0pt}]
            
            \node [rect node] (s11) {$1$};
            \node [rect node] (s12) [right=0cm of s11] {$4$};
            \node [rect node] (s13) [right=0cm of s12] {$2$};
            \node [rect node] (s14) [right=0cm of s13] {$3$};
            
            \node [rect node] (s21) [below=0cm of s11] {$3$};
            \node [rect node] (s22) [right=0cm of s21] {$2$};
            \node [rect node] (s23) [right=0cm of s22] {$4$};
            \node [rect node] (s24) [right=0cm of s23] {$1$};
            
            \node [rect node] (s31) [below=0cm of s21] {$4$};
            \node [rect node] (s32) [right=0cm of s31] {$1$};
            \node [rect node] (s33) [right=0cm of s32] {$3$};
            \node [rect node] (s34) [right=0cm of s33] {$2$};
            
            \node [rect node] (s41) [below=0cm of s31] {$2$};
            \node [rect node] (s42) [right=0cm of s41] {$3$};
            \node [rect node] (s43) [right=0cm of s42] {$1$};
            \node [rect node] (s44) [right=0cm of s43] {$4$};
            
        \end{tikzpicture}
    \end{minipage}%
    \begin{minipage}{.3\textwidth}
        \begin{tikzpicture}
            [->,auto,node distance=1.5cm,
            empt node/.style={font=\sffamily,inner sep=0pt,minimum size=0pt},
            rect node/.style={rectangle,draw,font=\sffamily,minimum size=0.8cm, inner sep=0pt, outer sep=0pt}]
            
            \node [rect node] (s11) {$1,1$};
            \node [rect node] (s12) [right=0cm of s11] {$3,4$};
            \node [rect node] (s13) [right=0cm of s12] {$4,2$};
            \node [rect node] (s14) [right=0cm of s13] {$2,3$};
            
            \node [rect node] (s21) [below=0cm of s11] {$4,3$};
            \node [rect node] (s22) [right=0cm of s21] {$2,2$};
            \node [rect node] (s23) [right=0cm of s22] {$1,4$};
            \node [rect node] (s24) [right=0cm of s23] {$3,1$};
            
            \node [rect node] (s31) [below=0cm of s21] {$2,4$};
            \node [rect node] (s32) [right=0cm of s31] {$4,1$};
            \node [rect node] (s33) [right=0cm of s32] {$3,3$};
            \node [rect node] (s34) [right=0cm of s33] {$1,2$};
            
            \node [rect node] (s41) [below=0cm of s31] {$3,2$};
            \node [rect node] (s42) [right=0cm of s41] {$1,3$};
            \node [rect node] (s43) [right=0cm of s42] {$2,1$};
            \node [rect node] (s44) [right=0cm of s43] {$4,4$};
            
        \end{tikzpicture}
    \end{minipage}%
    \caption{Orthogonal Latin squares of order $N=4$, and their superposition.}
    \label{fig:ols-ex}
\end{figure}
Contrarily to the case of single Latin squares, it is not possible to construct a family of MOLS for any possible order. A curious result is that the only two orders for which there are no orthogonal Latin squares are $N=2$ and $N=6$~\cite{stinson-cd}.

A second type of combinatorial designs that we consider in this paper are \emph{orthogonal arrays}, closely related to orthogonal Latin squares:
\begin{definition}
    \label{def:oa}
    A $(N, k, s, t)$ \emph{orthogonal array} ($(N, k, s, t)$--OA) is a $N \times k$ matrix with entries from a finite set $X$ of $s$ symbols such that, for any subset of $t$ columns, every $t$--uple of symbols occurs exactly $\lambda = N/s^t$ times.
\end{definition}
The parameter $t$ is also called the \emph{strength} of the OA. When $t=2$ and $\lambda=1$, the resulting orthogonal array is a $N^2\times k$ matrix in which every pair of columns contains all ordered pairs of symbols from $X$. In this case, the orthogonal array is simply denoted as $OA(k,v)$, and it is equivalent to a set of $(k-2)$--MOLS. We refer the reader to~\cite{stinson-cd,hedayat12} for a proof of this equivalence.

MOLS families and orthogonal arrays have several applications in cryptography and coding theory. Perhaps the best known example is the construction of threshold secret sharing schemes. A \emph{secret sharing scheme} (SSS) is a protocol that enables a dealer to share a secret value $S$ among a set of players $P_1, \ldots, P_n$ by distributing to them shares of $S$. The sharing operation is done in such a way that, at a later stage, only certain \emph{authorized subsets} can reconstruct the secret $S$ by pooling together their respective shares. In particular, if all non-authorized subsets do not gain any information on the value of the secret by combining their shares, the scheme is called \emph{perfect}. In a $(t,n)$ \emph{threshold} SSS, the authorized subsets are all those with at least $t$ participants~\cite{shamir79}. It can be shown that a $(t,n)$ perfect threshold secret sharing scheme is equivalent to an OA with strength $t$, $\lambda=1$ and $k=n+1$ columns~\cite{stinson-cd}. Accordingly, $(2,n)$ threshold schemes are equivalent to families of $n$ MOLS.

The link between MOLS and OA with coding theory is even deeper. As discussed for example in~\cite{hedayat12}, \emph{linear OA} (i.e. OA whose rows form a vector subspace) correspond to the \emph{duals} of linear codes. In particular, the strength $t$ of an OA is always strictly less than the dual distance $d^\bot$ of the corresponding linear code. Moreover, OA with index $\lambda=1$ are equivalent to \emph{Maximum Distance Separable} (MDS) codes, which satisfy the Singleton bound with equality. Due to the optimal distance of their codewords, MDS codes are also employed in symmetric cryptography when designing the so-called linear layer of a block cipher, since they achieve good diffusion properties~\cite{daemen20}.

\section{CA as Algebraic Systems}
\label{sec:ca-alg}

In this section, we describe the basic method to interpret any one-dimensional CA as an algebraic system, namely the \emph{block transformation}, and we survey the works that initially used this technique to study the dynamical properties of CA. Next, we show a related aspect of the block transformation, i.e. the preimage computation algorithm for permutive CA, and how it has been employed for cryptographic applications. Finally, we show how the algebraic structure induced by the block transformation of a bipermutive CA generates a Latin square.

\subsection{The Block Transformation}
\label{subsec:block}

Following Definition~\ref{def:nbca}, a CA is a vectorial function $F: \Sigma^{n} \to \Sigma^{n+d-1}$ where each output coordinate is determined by the application of the local rule $f$ on the corresponding neighborhood of diameter $d$. The \emph{block transformation} allows one to re-define any CA of diameter $d$ as a CA of diameter $2$, where each cell looks only at itself and its right neighbor to compute the next state.

\begin{figure}[t]
    \centering
    \begin{tikzpicture}
        [->,auto,node distance=1.5cm, empt node/.style={font=\sffamily,inner
            sep=0pt}, rect
        node/.style={rectangle,draw,font=\bfseries,minimum size=0.5cm, inner
            sep=0pt, outer sep=0pt}]
        
        \node [empt node] (c)   {};
        \node [rect node] (c1) [right=0.1cm of c] {$1$};
        \node [rect node] (c2) [right=0cm of c1] {$1$};
        \node [rect node] (c3) [right=0cm of c2] {$0$};
        \node [rect node] (c4) [right=0cm of c3] {$1$};
        
        \node [rect node] (c6) [right=3cm of c4, minimum width = 2cm, very thick] {$1101$};
        
        \node [empt node] (f1) [above=0.5cm of c2.east] {\Large $\Downarrow$ \normalsize $f$};
        
        \node [rect node] (p2) [above=0.85cm of c1] {$0$};
        \node [rect node] (p1) [left=0cm of p2] {$1$};
        \node [rect node] (p0) [left=0cm of p1] {$0$};
        \node [rect node] (p3) [right=0cm of p2] {$0$};
        \node [rect node] (p4) [right=0cm of p3] {$1$};
        \node [rect node] (p5) [right=0cm of p4] {$0$};
        \node [rect node] (p6) [right=0cm of p5] {$1$};
        \node [rect node] (p9) [right=0cm of p6] {$1$};
        
        \node [rect node] (p10) [right=1cm of p9, minimum width = 2cm, very thick] {$0100$};
        \node [rect node] (p11) [right=0cm of p10, minimum width = 2cm, very thick] {$1011$};
        
        \node [empt node] (f3) [right=1.8cm of f1] {\LARGE $\Rightarrow$};
        
        \node [empt node] (f2) [right=4.3cm of f1] {\Large $\Downarrow$ \normalsize $f'$};
        
        \node [empt node] (p7) [below=0.2cm of p1] {};
        \node [empt node] (p8) [right=0.07cm of p7] {};
        \node [empt node] (p12) [above=0.5cm of p1.east] {};
        \node [empt node] (p13) [above=0.5cm of p5.east] {};
        \node [empt node] (p14) [above=0.3cm of p13] {\phantom{M}};
        
        \draw[-,very thick] (p0.north west) -- (p9.north east);
        \draw[-,very thick] (p0.south west) -- (p9.south east);
        \draw[-,very thick] (p0.north west) -- (p0.south west);
        \draw[-,very thick] (p9.north east) -- (p9.south east);
        \draw[-,very thick] (p3.north east) -- (p3.south east);
        
        \draw[-,very thick] (c1.north west) -- (c4.north east);
        \draw[-,very thick] (c1.south west) -- (c4.south east);
        \draw[-,very thick] (c1.north west) -- (c1.south west);
        \draw[-,very thick] (c4.north east) -- (c4.south east);
        
    \end{tikzpicture}
    \caption{Example of block transformation for a CA of diameter $d=5$. The original local rule is defined as $f(x_1,x_2,x_3,x_4,x_5) = x_1 \oplus x_3 \oplus x_5$.}
    \label{fig:block}
\end{figure}
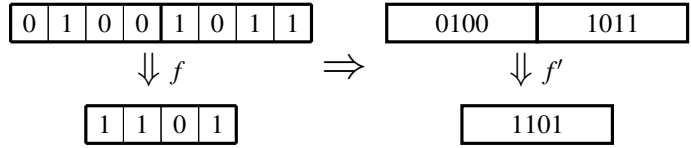

The idea is to group the cells in blocks of length $d-1$, redefining the CA's alphabet as $\hat{\Sigma} = \Sigma^{d-1}$. The new local rule $f': \hat{\Sigma}^2 \to \hat{\Sigma}$ thus maps two $(d-1)$-cell blocks to a single $(d-1)$-cell block. This can be simulated with the original local rule $f$ by just defining $f'$ as the NBCA $G: \Sigma^{2(d-1)} \to \Sigma^{d-1}$, where $f$ is applied to an input vector of $2(d-1)$ cells, so that the final output vector is composed of $d-1$ cells. Figure~\ref{fig:block} depicts an example of the block transformation on a CA of diameter $d=5$. Assuming that the original alphabet is binary, the local rule is a Boolean function of the form $f: \F_2^5 \to \F_2$. If we arrange the cells four by four, the new alphabet is $\hat{\Sigma} =  \F_2^4$, i.e. 4-tuples of bits. The new local rule $f'$ will thus have the form $f': (\F_2^4)^2 \to \F_2^4$.

An advantage of the block transformation is that the resulting local rule $f'$ of diameter $d=2$ can be interpreted as a \emph{binary operation} $*: \hat{\Sigma} \times \hat{\Sigma} \to \hat{\Sigma}$. Hence, one can consider a CA essentially as an algebraic structure $\langle \hat{\Sigma}, * \rangle$ over the state alphabet. The earlier literature that studied CA as algebraic systems focused on the relationship between the dynamic behavior of an infinite CA and the properties satisfied by the underlying algebraic structure $\langle \hat{\Sigma}, * \rangle$, such as being a group or a monoid.

The block transformation was introduced by Pedersen in~\cite{p-complexsyst-1992}\footnote{Strictly speaking, the term ``block transformation'' actually comes from a paper by Moore and Drisko~\cite{md-complexsyst-1996}, although Pedersen already used it without giving it a name. The same technique was sketched even earlier by Albert and Cul{\'{\i}}k~\cite{albert87} and Smith~\cite{smith71}, as well as by Hedlund in his seminal work on CA~\cite{hedlund69}.}, where he used it to characterize infinite CA with ultimately periodic behavior to local rules that are subvarieties of groupoids. In this context, a ``groupoid'' is simply an algebraic structure $\langle \Sigma, * \rangle$ whose only property is closure, i.e., the result of the operation $*$ is still an element of the set to which the two operands belong to\footnote{The terminology here is not standard. In abstract algebra, groupoids are usually called \emph{magmas}.}. The author then investigated the relationship between different varieties of these groupoids.

Eloranta~\cite{e-nonlin-1993} considered the block transformation defined by partially permutive local rules (i.e., rules that are permutive only on subsets of the state alphabet) to explain kink-like structures exhibited by the dynamic evolution of certain CA. To the best of our knowledge, this is the first work where the block transformation of a permutive CA is related to a quasigroup.

The algebraic perspective on CA has been further brought forward in several works by Moore and coworkers. Moore and Drisko~\cite{md-complexsyst-1996} investigated which algebraic structures give rise to efficiently predictable CA. In particular, they showed that if the binary operation defined by the block transformation satisfies one of four properties (associativity with identity, inverse property loop, anticommutativity with identity and commutativity), then the original CA local rule depends only on its leftmost and rightmost cells. Moore~\cite{m-physicad-1997} showed that for a variety of algebraic structures such as quasigroups and Steiner systems, there exists an efficient algorithm to predict the dynamical evolution of the corresponding CA. For this reason, he termed such CA as \emph{quasilinear}, since they obey a law analogous to the superposition principle of linear CA. Moore and Boykett~\cite{moore97} considered the problem of \emph{commuting} cellular automata, that is, under which conditions one can apply the global rules of two CA in any order and obtain the same output configuration. The main finding of that work is that linear permutive CA cannot commute with nonlinear CA. Finally, Moore~\cite{m-physicad-1998} proved that nonlinear CA whose block transformation yields a solvable group can be decomposed into the quasidirect product of linear CA, making their prediction by parallel circuits more efficient.

\subsection{Preimage Computation for Permutive CA}
\label{subsec:preim}
In this section we discuss a research thread that developed alongside the perspective of CA as algebraic systems, which has been specifically adopted for the design of cryptographic primitives.

Since Hedlund's work~\cite{hedlund69}, it is well-known that permutive CA are surjective. The reasoning goes as follows: suppose that we have a finite configuration $y \in \Sigma^{n-d+1}$ and we want to determine one of its preimages under the action of a permutive CA $F: \Sigma^n \to \Sigma^{n-d+1}$. Without loss of generality, let us assume that the local rule $f: \Sigma^d \to \Sigma$ is rightmost permutive, and fix the leftmost $(d-1)$-cell block of a preimage $x \in \Sigma^{n}$ to an arbitrary value $x_{1\ldots d-1} = (x_1,\ldots, x_{d-1}) \in \Sigma^{d-1}$. Then, by rightmost permutivity we know that $x_{1\ldots d-1}$ determines a permutation $\pi: \Sigma \to \Sigma$ between $x_{d}$ and the value of the output cell $y_1$. Therefore, we can compute $x_{d}$ by applying the inverse permutation $\pi^{-1}$ to $y_1$. After that, the $(d-1)$-cell block $x_{2\ldots d} = (x_2, \ldots, x_{d})$ is fully determined, so we can obtain $x_{d+1}$ as $\pi^{-1}(y_2)$; and so on, until the preimage is complete. A similar argument applies if the rule is leftmost permutive, by initializing the rightmost $(d-1)$-cell block of the preimage and then completing it to the left. If the rule is bipermutive the initial block of $d-1$ cells can be put in any position of the preimage, which is then completed by ``expanding'' the block in both directions.

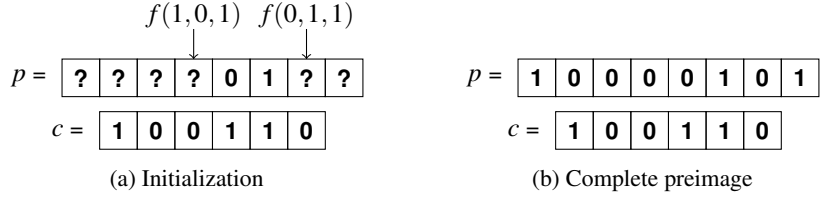
\begin{figure}[t]
    \centering
    \begin{subfigure}{.5\textwidth}
        \centering
        \begin{tikzpicture}
            [->,auto,node distance=1.5cm,
            empt node/.style={font=\sffamily,inner sep=0pt,minimum size=0pt},
            rect node/.style={rectangle,draw,font=\sffamily\bfseries,minimum size=0.5cm, inner sep=0pt, outer sep=0pt}]
            
            \node [empt node] (c)   {$c$ =};
            \node [rect node] (c1) [right=0.2cm of c] {1};
            \node [rect node] (c2) [right=0cm of c1] {0};
            \node [rect node] (c3) [right=0cm of c2] {0};
            \node [rect node] (c4) [right=0cm of c3] {1};
            \node [rect node] (c5) [right=0cm of c4] {1};
            \node [rect node] (c6) [right=0cm of c5] {0};
            
            \node [rect node] (p2) [above=0.2cm of c1] {?};
            \node [rect node] (p1) [left=0cm of p2] {?};
            \node [empt node] (p)  [left=0.2cm of p1] {$p$ =};
            \node [rect node] (p3) [right=0cm of p2] {?};
            \node [rect node] (p4) [right=0cm of p3] {?};
            \node [rect node] (p5) [right=0cm of p4] {0};
            \node [rect node] (p6) [right=0cm of p5] {1};
            \node [rect node] (p7) [right=0cm of p6] {?};
            \node [rect node] (p8) [right=0cm of p7] {?};
            
            \node [empt node] (f1) [above=0.4cm of p4] {$f(1,0,1)$};
            \node [empt node] (f2) [above=0.4cm of p7] {$f(0,1,1)$};
            
            \draw[->] (f1) -- (p4);
            \draw[->] (f2) -- (p7);
            
        \end{tikzpicture}
        \caption{Initialization}
    \end{subfigure}%
    \begin{subfigure}{.5\textwidth}
        \centering
        \begin{tikzpicture}
            [->,auto,node distance=1.5cm,
            empt node/.style={font=\sffamily,inner sep=0pt,minimum size=0.3cm},
            rect node/.style={rectangle,draw,font=\sffamily\bfseries,minimum size=0.5cm, inner sep=0pt, outer sep=0pt}]
            
            \node [empt node] (c)   {$c$ =};
            \node [rect node] (c1) [right=0.2cm of c] {1};
            \node [rect node] (c2) [right=0cm of c1] {0};
            \node [rect node] (c3) [right=0cm of c2] {0};
            \node [rect node] (c4) [right=0cm of c3] {1};
            \node [rect node] (c5) [right=0cm of c4] {1};
            \node [rect node] (c6) [right=0cm of c5] {0};
            
            \node [rect node] (p2) [above=0.2cm of c1] {0};
            \node [rect node] (p1) [left=0cm of p2] {1};
            \node [empt node] (p)  [left=0.2cm of p1] {$p$ =};
            \node [rect node] (p3) [right=0cm of p2] {0};
            \node [rect node] (p4) [right=0cm of p3] {0};
            \node [rect node] (p5) [right=0cm of p4] {0};
            \node [rect node] (p6) [right=0cm of p5] {1};
            \node [rect node] (p7) [right=0cm of p6] {0};
            \node [rect node] (p8) [right=0cm of p7] {1};
            
            \node [empt node] (f1) [above=0.4cm of p4] {};
        \end{tikzpicture}
        \caption{Complete preimage}
    \end{subfigure}
    \caption{Preimage computation for $c=(1,0,0,1,1,0) \in \mathbb{F}_2^6$ using rule 150.}
    \label{fig:bip-inv-ex}
\end{figure}

The procedure above is also called the \emph{preimage computation} or \emph{preimage reconstruction} algorithm by some authors~\cite{macedo08,ml-acri-2014}. This procedure also has an elegant description in terms of de Bruijn graphs: the basic idea is to find the path on the edges labeled by the output configuration, and then return the corresponding path on the vertices (merged together with the fusion operator) as the preimage~\cite{sutner91,mariot15}. As an example, Figure~\ref{fig:bip-inv-ex} depicts the preimage reconstruction process for a binary NBCA $F: \F_2^8 \to \F_2^6$ equipped with the bipermutive rule 150.	
In the binary case, the inverse permutation $\pi^{-1}$ used to determine the value of the leftmost (respectively, rightmost) input cell is simply the XOR between the output value, the rightmost (respectively, leftmost) input cell and the generating function $g$ computed on the central $d-2$ input cells. As a matter of fact, one has:
\begin{equation}
    \label{eq:inv-bip}
    y = x_1 \oplus g(x_2,\ldots,x_{d-1}) \oplus x_d \ \Leftrightarrow \ x_1 = y \oplus  g(x_2,\ldots,x_{d-1}) \oplus x_d \enspace ,
\end{equation}
or equivalently,
\begin{equation}
    \label{eq:inv-bip-2}
    y = x_1 \oplus g(x_2,\ldots,x_{d-1}) \oplus x_d \ \Leftrightarrow \ x_d = x_1 \oplus  g(x_2,\ldots,x_{d-1}) \oplus y \enspace .
\end{equation}

Remark that, in the no-boundary setting, the preimage computation algorithm can be iterated indefinitely, obtaining at each step a slightly larger preimage that has $d-1$ additional cells. This observation has been used by various researchers to devise cryptographic applications based on permutive CA, the first one being Gutowitz~\cite{gutowitz93}. There, the author proposed to iterate the preimage computation algorithm with a permutive CA to design the diffusion phase of a block cipher. Oliveira et al.~\cite{oliveira04} refined Gutowitz's idea by considering bipermutive rules (there called ``bi-directional toggle rules''). The underlying argument for using bipermutive rules instead of just leftmost or rightmost ones is that differences propagate in both directions, making differential cryptanalysis more difficult. The main problem of the proposal, however, remained ciphertext expansion as in Gutowitz's original design, which is a direct consequence of the iterated application of the preimage reconstruction algorithm. Later, Mac{\^{e}}do et al.~\cite{macedo08} attempted to address this issue by forcing reversibility with a periodic-boundary CA.

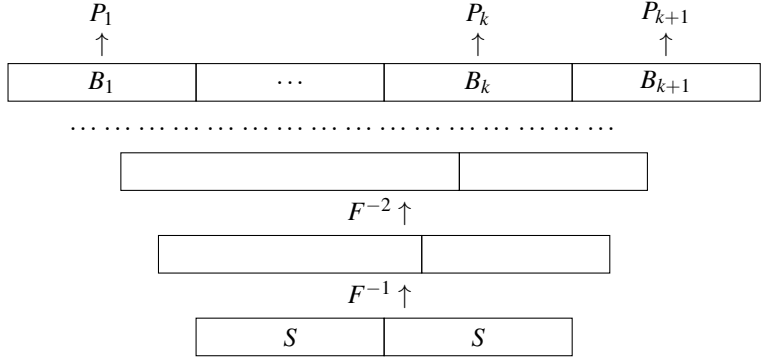
\begin{figure}[t]
    \centering
    \begin{tikzpicture}
        [->,auto,node distance=1.5cm, empt node/.style={font=\sffamily,inner
            sep=0pt,outer sep=0pt}, rect
        node/.style={rectangle,draw,font=\sffamily\bfseries,minimum
            width=1cm,minimum height=0.5cm, inner sep=0pt, outer sep=0pt}, white
        node/.style={rectangle,draw,font=\sffamily\bfseries,minimum width=1cm,minimum
            height=0.5cm, inner sep=0pt, outer sep=0pt}, null
        node/.style={rectangle,draw=white,font=\sffamily\bfseries,minimum width=1cm,minimum
            height=0.5cm, inner sep=0pt, outer sep=0pt}]
        
        \node [white node, minimum width=2.5cm] (s) {$S$};
        \node [rect node, minimum width=3.5cm] (p1) [above=0.6cm of s] {};
        \node [empt node] (e1) [above=0.3cm of s] {};
        \node [empt node] (e2) [right=0.75cm of e1] {$F^{-1}\uparrow$};
        \node [empt node] (e3) [above=0.3cm of p1] {};
        \node [empt node] (e4) [right=0.75cm of e3] {$F^{-2}\uparrow$};
        \node [rect node, minimum width=4.5cm] (p2) [above=0.6cm of p1] {};
        \node [empt node] (e5) [above=2.5cm of s.east] {};
        \node [empt node] (dots1) [above=2.5cm of s] {$\ldots$ $\ldots$ $\ldots$
            $\ldots$ $\ldots$ $\ldots$ $\ldots$ $\ldots$ $\ldots$ $\ldots$ $\ldots$
            $\ldots$ $\ldots$};
        \node [rect node, minimum width=2.5cm] (p32) [above=0.35cm of dots1] {$\ldots$};
        \node [rect node, minimum width=2.5cm] (p31) [left=0cm of p32] {$B_1$};
        \node [rect node, minimum width=2.5cm] (p33) [right=0cm of p32] {$B_k$};
        
        \node [white node, minimum width=2.5cm] (s1) [right=0cm of s] {$S$};
        \node [white node, minimum width=2.5cm] (p12) [right=0cm of p1] {};
        \node [white node, minimum width=2.5cm] (p22) [right=0cm of p2] {};
        \node [empt node] (dots2) [right=0.1cm of dots1] {$\ldots$ $\ldots$ $\ldots$};
        \node [white node, minimum width=2.5cm] (p34) [right=0cm of p33] {$B_{k+1}$};
        
        \node [empt node] (arr1) [above=0.1cm of p31] {$\uparrow$};
        \node [empt node] (arrk) [above=0.1cm of p33] {$\uparrow$};
        \node [empt node] (arrk1) [above=0.1cm of p34] {$\uparrow$};
        \node [empt node] (pl1) [above=0.1cm of arr1] {$P_1$};
        \node [empt node] (plk) [above=0.1cm of arrk] {$P_k$};
        \node [empt node] (plk1) [above=0.1cm of arrk1] {$P_{k+1}$};
        
        \node [empt node] (aw12) [above right=0cm of p1] {};
        \node [empt node] (aw22) [above right=0cm of p2] {};
        \node [empt node] (aw32) [above right=0cm of p33] {};
        
    \end{tikzpicture}
    \caption{Setup phase of the SSS scheme from~\cite{ml-acri-2014} with two copies of the secret $S$.}
    \label{fig:setup-phase}
\end{figure}

More recently, Mariot and Leporati~\cite{ml-acri-2014} introduced a CA-based secret sharing scheme based on the iterated preimage computation algorithm. The secret, denoted as $S$, is represented as a finite configuration of length $m$ in a NBCA equipped with a bipermutive local rule. The dealer executes the preimage construction algorithm until a preimage of length $k \cdot m$ is achieved, where $k$ is the number of players.	This preimage is then divided into $k$ blocks, each of size $m$, and sequentially distributed to the players. To reconstruct the secret, players must correctly order and combine their blocks, subsequently evolving the CA forward until the original secret is obtained. In this way, the resulting scheme is a $(k,k)$-threshold SSS, where all $k$ players are required to pool their shares to recover the secret. A straightforward modification extends this scheme to accommodate $k + 1$ players. This involves appending a \emph{copy} of the secret to the right, resulting in a final preimage with $k + 1$ blocks of size $m$. Both sets of players, namely $P_1, \dots,P_k$ and $P_2,\dots,P_{k+1}$, can then recover the secret using the same procedure: by combining their respective shares, the dynamic evolution of the resulting partial preimage collapses on one of the two copies of the secret. Figure~\ref{fig:setup-phase} displays the setup phase of the SSS in this specific case. 

This procedure can be generalized to establish a $(k, n)$-threshold scheme by concatenating $k$ copies of the secret. Consequently, the scheme exhibits a \emph{sequential} threshold access structure, where all minimal authorized subsets take the form $\{P_i,\dots,P_{i+k-1}\}$. This sequentiality feature is shared by other CA-based SSS, such as the one proposed by del Rey et al. in~\cite{rey05}. However, a notable departure from the approach in \cite{rey05} lies in the fact that, in the latter, shares must adhere to a \emph{temporal adjacency} constraint, being successive configurations of a higher-order CA. On the other hand, in the scheme proposed in~\cite{ml-acri-2014} the shares are \emph{spatially adjacent}, since they constitute blocks of an NBCA preimage.

\begin{figure}[t]
    \centering
    \begin{tikzpicture}
        [->,auto,node distance=1.5cm,
        empt node/.style={font=\sffamily,inner sep=0pt,outer sep=0pt},
        rect node/.style={rectangle,draw,font=\sffamily\bfseries,minimum width=1cm,minimum height=0.6cm, inner sep=0pt, outer sep=0pt},
        white node/.style={rectangle,draw=white,fill=white,minimum width=0.2cm,minimum height=0.6cm, inner sep=0pt, outer sep=0pt}]
        
        \node [rect node, minimum width=2.5cm] (s1) {$S$};
        \node [rect node, minimum width=0.75cm] (s0) [left=0cm of s1] {};
        \node [empt node] (a1) {};
        \node [white node] (c1) [left=0cm of s0] {};
        \node [empt node] (d1) [left=1.75cm of a1] {$\ldots$};
        \node [rect node, minimum width=2.5cm] (s2) [right=0cm of s1] {$\ldots$};
        \node [rect node, minimum width=2.5cm] (s3) [right=0cm of s2] {$S$};
        \node [rect node, minimum width=2.5cm] (s4) [right=0cm of s3] {$S$};
        \node [rect node, minimum width=0.75cm] (s5) [right=0cm of s4] {};
        \node [white node] (c3) [right=0cm of s5] {};
        \node [empt node] (a3) [right=0cm of s4] {};
        \node [empt node] (d3) [right=0.5cm of a3] {$\ldots$};
        \node [empt node] (e2) [left=0cm of s1] {};
        \node [empt node] (e1) [above=0.8cm of s1] {};
        \node [rect node] (w1) [left=0.75cm of e1] {$w_1$};
        \node [rect node, minimum width=0.25cm] (w0) [left=0cm of w1] {};
        \node [empt node] (a2) [left=0.75cm of e1] {};
        \node [white node] (c2) [left=0cm of w0] {};
        \node [empt node] (d2) [left=1cm of a2] {$\ldots$};
        \node [rect node, minimum width=1.5cm] (b1) [right=0cm of w1] {$v_1$};
        \node [rect node, minimum width=1cm] (w2) [right=0cm of b1] {$w_2$};
        \node [rect node, minimum width=1.5cm] (b2) [right=0cm of w2] {$\ldots$};
        \node [rect node] (w3) [right=0cm of b2] {$w_{h-1}$};
        \node [rect node, minimum width=1.5cm] (b3) [right=0cm of w3] {$v_{h-1}$};
        \node [rect node] (w4) [right=0cm of b3] {$w_1$};
        \node [rect node, minimum width=1.5cm] (b4) [right=0cm of w4] {$v_1$};
        \node [rect node] (w5) [right=0cm of b4] {$w_2$};
        \node [rect node, minimum width=0.25cm] (b5) [right=0cm of w5] {};
        \node [empt node] (a4) [right=0cm of w5] {};
        \node [white node] (c4) [right=0.25cm of w5] {};
        \node [empt node] (d4) [right=0cm of a4] {$\ldots$};
        
        \draw [-, decorate, decoration={brace,amplitude=5pt,mirror,raise=0.4cm}]
        (e2) -- (s4) node [midway,yshift=-1.1cm] {$h \le |\Sigma|^{d-1}$ copies of
            $s$};
        
        \draw[->, thick, shorten >=2pt,shorten <=2pt,>=stealth] (w1) edge[bend left] node (f1) [above]{} (w2.north);
        \draw[->, dashed, thick, shorten >=2pt,shorten <=2pt,>=stealth] (w2) edge[bend left] node (f2) [above]{}  (w3.north);
        \draw[->, thick, shorten >=2pt,shorten <=2pt,>=stealth] (w3) edge[bend left] node (f3) [above]{} (w4.north);
        \draw[->, thick, shorten >=2pt,shorten <=2pt,>=stealth] (w4) edge[bend left] node (f4) [above]{} (w5.north);
        
    \end{tikzpicture}
    \caption{Repetition of shares in the SSS proposed in~\cite{ml-acri-2014}.}
    \label{fig:per-preim}
\end{figure}
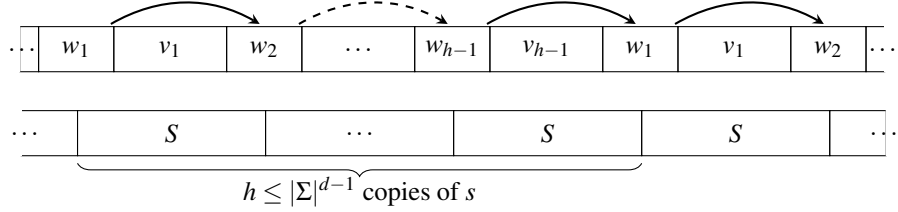

Furthermore, the authors of~\cite{ml-acri-2014} remark that, as the preimage of a configuration is uniquely determined by a block of $d-1$ cells, the shares must eventually repeat after at most $|\Sigma|^{d-1}$ juxtaposed copies of the secret. Hence, the resulting access structure of the scheme is both sequential and \emph{cyclic}. Figure~\ref{fig:per-preim} represents this effect for a single iteration of the preimage computation algorithm. The $d-1$-cell blocks that determine the subsequent block of $m$ cells are denoted as $w_1,\ldots,w_{h-1}$.

Consequently, determining the maximum number of players allowed in the CA-based SSS of~\cite{ml-acri-2014} is equivalent to studying the periods of preimages of \emph{spatially periodic configurations} (SPC) in bipermutive CA. Indeed, a well-known fact for infinite CA is that any preimage $x \in F^{-1}(y)$ of a SPC $y \in \Sigma^{\Z}$ is also a SPC, whose period is a multiple of the period of $y$~\cite{hedlund69}. Following the SSS-inspired motivation, the authors of~\cite{mariot15} investigated more in detail the periods of preimages in SPC, first by providing some upper bounds for the case of general bipermutive CA. Then, restricting the attention to the subclass of linear bipermutive CA over a finite field $\F_q$, the authors provided an exact characterization of the periods of preimages, leveraging on the theory of \emph{Linear Recurring Sequences} (LRS). Indeed, the preimage reconstruction algorithm for a linear bipermutive CA can be synthesized by a concatenation of \emph{Linear Feedback Shift Registers} (LFSR), and thus the period of the preimage can be deduced by the minimal polynomials of the two LFSRs.

\subsection{Latin Squares from Bipermutive CA}
\label{subsec:ls-ca}

We conclude this section by showing how the two research tracks presented so far---CA as algebraic systems and preimage computation for permutive CA--- are actually related to one another, and provide the basis for the combinatorial designs perspective that we develop in the next section.

From now on, we focus on the case of an NBCA $F: \Sigma^{2(d-1)} \to \Sigma^{d-1}$ defined by a bipermutive local rule $f: \Sigma^d \to \Sigma$ of diameter $d$. In this way, we can directly consider the CA global rule as an algebraic structure $\langle \hat{\Sigma}, * \rangle$, where $\hat{\Sigma} = \Sigma^{d-1}$ and $* \equiv F$. In general, the \emph{Cayley table} of a finite algebraic structure $\langle S, * \rangle$ with $N = |S|$ is the $N \times N$ matrix $C_*$ where the rows and columns are indexed by the elements of $S$, and $C_*(x, y) = x * y$ for all pairs of row and column coordinates $x,y \in S$. For the CA case $\langle \hat{\Sigma}, F \rangle$, denoting $N = |\hat{\Sigma}| = |\Sigma|^{d-1}$, the corresponding Cayley table $C_F$ is defined by using the leftmost $(d-1)$-cell block as the row coordinate, the rightmost $(d-1)$-cell block as the column coordinate, and the output configuration of the CA (itself a $(d-1)$-cell block) as the entry to be inserted at those coordinates.

Formally, let us define a total order $\le$ over $\Sigma^{d-1}$, and let $\phi: \Sigma^{d-1} \rightarrow [N]$ be a
monotone bijective mapping between $\Sigma^{d-1}$ and $[N] = \{1,\ldots,q^{d-1}\}$,	where the total order on $[N]$ is inherited from the usual order of natural numbers. We denote by $\psi = \phi^{-1}$, i.e. the inverse of $\phi$. We can now give the formal definition of the Cayley table associated to a CA:
\begin{definition}
    \label{def:cayley}
    Let $\Sigma$ be an alphabet of $q$ symbols and $d \in \N$, with $N = q^{d-1}$. Further, suppose that $F: \Sigma^{2(d-1)} \rightarrow \Sigma^{d-1}$ is a CA defined by the local rule $f: \Sigma^d \rightarrow \Sigma$. The \emph{Cayley table} associated to $F$ is the $N \times N$ matrix $C_F$ with entries from $[N]$ such that
    \begin{equation}
        \label{eq:sq-ca}
        C_{F}(i,j) = \phi(F(\psi(i)||\psi(j))) \enspace ,
    \end{equation}
    for all $1 \le i,j \le N$, where $\psi(i)||\psi(j) \in \Sigma^{2(d-1)}$ denotes the \emph{concatenation} of  $\psi(i),\psi(j) \in \Sigma^{d-1}$.
\end{definition}

The next result, originally proved by Eloranta~\cite{e-nonlin-1993} and independently re-discovered by Mariot et al.~\cite{mgfl-desi-2020}, shows the connection between bipermutive CA, Latin squares and quasigroups:
\begin{lemma}[\cite{e-nonlin-1993,mgfl-desi-2020}]
\label{lm:lat-sq-bip-ca}
    The Cayley table $C_F$ of a CA $F: \Sigma^{2(d-1)} \rightarrow \Sigma^{d-1}$ defined by a bipermutive local rule
    $f:\Sigma^{d}\rightarrow \Sigma$ is a Latin square of order $N=q^{d-1}$, where $q = |\Sigma|$. Equivalently, the algebraic structure $\langle \hat{\Sigma}, F \rangle$ is a quasigroup of order $N$.
\end{lemma}
The proof of this lemma is an application of the preimage computation algorithm in both directions: by fixing the leftmost block of the CA, there exists a permutation between the rightmost block and the output block. This means that by fixing any row of the Cayley table, we see all numbers from $1$ to $N$ exactly once. A symmetric argument holds by fixing the rightmost block of the CA (which means fixing a column of the table).

Figure~\ref{fig:r150-sq} illustrates the construction of the Cayley table associated to the CA $F: \F_2^4 \rightarrow \F_2^2$ defined by rule $150$. The mapping $\phi$ here is $\phi(00) \mapsto 1$, $\phi(10) \mapsto 2$, $\phi(01) \mapsto 3$ and
$\phi(11) \mapsto 4$. Indeed, one can see that each number from $1$ to $4$ occurs exactly once in each row and column of the table, making it a Latin square of order $4$.

\begin{figure}[t]
    \centering
    \begin{subfigure}{.5\textwidth}
        \centering
        \begin{tikzpicture}
            [->,auto,node distance=1.5cm,
            empt node/.style={font=\sffamily,inner sep=0pt,minimum size=0pt},
            rect node/.style={rectangle,draw,font=\sffamily,minimum size=0.3cm, inner sep=0pt, outer sep=0pt}]
            
            \node [empt node] (e1) {};
            \node [rect node] (i111) [right=0.5cm of e1] {$0$};
            \node [rect node] (i112) [right=0cm of i111] {$0$};
            \node [rect node] (i113) [right=0cm of i112] {$0$};
            \node [rect node] (i114) [right=0cm of i113] {$0$};
            \node [rect node] (i115) [below=0cm of i112] {$0$};
            \node [rect node] (i116) [right=0cm of i115] {$0$};
            
            \node [rect node] (i121) [right=0.3cm of i114] {$0$};
            \node [rect node] (i122) [right=0cm of i121] {$0$};
            \node [rect node] (i123) [right=0cm of i122] {$1$};
            \node [rect node] (i124) [right=0cm of i123] {$0$};
            \node [rect node] (i125) [below=0cm of i122] {$1$};
            \node [rect node] (i126) [right=0cm of i125] {$1$};
            
            \node [rect node] (i131) [right=0.3cm of i124] {$0$};
            \node [rect node] (i132) [right=0cm of i131] {$0$};
            \node [rect node] (i133) [right=0cm of i132] {$0$};
            \node [rect node] (i134) [right=0cm of i133] {$1$};
            \node [rect node] (i135) [below=0cm of i132] {$0$};
            \node [rect node] (i136) [right=0cm of i135] {$1$};
            
            \node [rect node] (i141) [right=0.3cm of i134] {$0$};
            \node [rect node] (i142) [right=0cm of i141] {$0$};
            \node [rect node] (i143) [right=0cm of i142] {$1$};
            \node [rect node] (i144) [right=0cm of i143] {$1$};
            \node [rect node] (i145) [below=0cm of i142] {$1$};
            \node [rect node] (i146) [right=0cm of i145] {$0$};
            
            \node [rect node] (i211) [below=0.6cm of i111] {$1$};
            \node [rect node] (i212) [right=0cm of i211] {$0$};
            \node [rect node] (i213) [right=0cm of i212] {$0$};
            \node [rect node] (i214) [right=0cm of i213] {$0$};
            \node [rect node] (i215) [below=0cm of i212] {$1$};
            \node [rect node] (i216) [right=0cm of i215] {$0$};
            
            \node [rect node] (i221) [right=0.3cm of i214] {$1$};
            \node [rect node] (i222) [right=0cm of i221] {$0$};
            \node [rect node] (i223) [right=0cm of i222] {$1$};
            \node [rect node] (i224) [right=0cm of i223] {$0$};
            \node [rect node] (i225) [below=0cm of i222] {$0$};
            \node [rect node] (i226) [right=0cm of i225] {$1$};
            
            \node [rect node] (i231) [right=0.3cm of i224] {$1$};
            \node [rect node] (i232) [right=0cm of i231] {$0$};
            \node [rect node] (i233) [right=0cm of i232] {$0$};
            \node [rect node] (i234) [right=0cm of i233] {$1$};
            \node [rect node] (i235) [below=0cm of i232] {$1$};
            \node [rect node] (i236) [right=0cm of i235] {$1$};
            
            \node [rect node] (i241) [right=0.3cm of i234] {$1$};
            \node [rect node] (i242) [right=0cm of i241] {$0$};
            \node [rect node] (i243) [right=0cm of i242] {$1$};
            \node [rect node] (i244) [right=0cm of i243] {$1$};
            \node [rect node] (i245) [below=0cm of i242] {$0$};
            \node [rect node] (i246) [right=0cm of i245] {$0$};
            
            \node [rect node] (i311) [below=0.6cm of i211] {$0$};
            \node [rect node] (i312) [right=0cm of i311] {$1$};
            \node [rect node] (i313) [right=0cm of i312] {$0$};
            \node [rect node] (i314) [right=0cm of i313] {$0$};
            \node [rect node] (i315) [below=0cm of i312] {$1$};
            \node [rect node] (i316) [right=0cm of i315] {$1$};
            
            \node [rect node] (i321) [right=0.3cm of i314] {$0$};
            \node [rect node] (i322) [right=0cm of i321] {$1$};
            \node [rect node] (i323) [right=0cm of i322] {$1$};
            \node [rect node] (i324) [right=0cm of i323] {$0$};
            \node [rect node] (i325) [below=0cm of i322] {$0$};
            \node [rect node] (i326) [right=0cm of i325] {$0$};
            
            \node [rect node] (i331) [right=0.3cm of i324] {$0$};
            \node [rect node] (i332) [right=0cm of i331] {$1$};
            \node [rect node] (i333) [right=0cm of i332] {$0$};
            \node [rect node] (i334) [right=0cm of i333] {$1$};
            \node [rect node] (i335) [below=0cm of i332] {$1$};
            \node [rect node] (i336) [right=0cm of i335] {$0$};
            
            \node [rect node] (i341) [right=0.3cm of i334] {$0$};
            \node [rect node] (i342) [right=0cm of i341] {$1$};
            \node [rect node] (i343) [right=0cm of i342] {$1$};
            \node [rect node] (i344) [right=0cm of i343] {$1$};
            \node [rect node] (i345) [below=0cm of i342] {$0$};
            \node [rect node] (i346) [right=0cm of i345] {$1$};
            
            \node [rect node] (i411) [below=0.6cm of i311] {$1$};
            \node [rect node] (i412) [right=0cm of i411] {$1$};
            \node [rect node] (i413) [right=0cm of i412] {$0$};
            \node [rect node] (i414) [right=0cm of i413] {$0$};
            \node [rect node] (i415) [below=0cm of i412] {$0$};
            \node [rect node] (i416) [right=0cm of i415] {$1$};
            
            \node [rect node] (i421) [right=0.3cm of i414] {$1$};
            \node [rect node] (i422) [right=0cm of i421] {$1$};
            \node [rect node] (i423) [right=0cm of i422] {$1$};
            \node [rect node] (i424) [right=0cm of i423] {$0$};
            \node [rect node] (i425) [below=0cm of i422] {$1$};
            \node [rect node] (i426) [right=0cm of i425] {$0$};
            
            \node [rect node] (i431) [right=0.3cm of i424] {$1$};
            \node [rect node] (i432) [right=0cm of i431] {$1$};
            \node [rect node] (i433) [right=0cm of i432] {$0$};
            \node [rect node] (i434) [right=0cm of i433] {$1$};
            \node [rect node] (i435) [below=0cm of i432] {$0$};
            \node [rect node] (i436) [right=0cm of i435] {$0$};
            
            \node [rect node] (i441) [right=0.3cm of i434] {$1$};
            \node [rect node] (i442) [right=0cm of i441] {$1$};
            \node [rect node] (i443) [right=0cm of i442] {$1$};
            \node [rect node] (i444) [right=0cm of i443] {$1$};
            \node [rect node] (i445) [below=0cm of i442] {$1$};
            \node [rect node] (i446) [right=0cm of i445] {$1$};
            
        \end{tikzpicture}
    \end{subfigure}%
    \begin{subfigure}{.5\textwidth}
        \centering
        \begin{tikzpicture}
            [->,auto,node distance=1.5cm,
            empt node/.style={font=\sffamily,inner sep=0pt,minimum size=0pt},
            rect node/.style={rectangle,draw,font=\sffamily,minimum size=0.8cm, inner sep=0pt, outer sep=0pt}]
            \large
            
            \node [rect node] (s11) {$1$};
            \node [rect node] (s12) [right=0cm of s11] {$4$};
            \node [rect node] (s13) [right=0cm of s12] {$3$};
            \node [rect node] (s14) [right=0cm of s13] {$2$};
            
            \node [rect node] (s21) [below=0cm of s11] {$2$};
            \node [rect node] (s22) [right=0cm of s21] {$3$};
            \node [rect node] (s23) [right=0cm of s22] {$4$};
            \node [rect node] (s24) [right=0cm of s23] {$1$};
            
            \node [rect node] (s31) [below=0cm of s21] {$4$};
            \node [rect node] (s32) [right=0cm of s31] {$1$};
            \node [rect node] (s33) [right=0cm of s32] {$2$};
            \node [rect node] (s34) [right=0cm of s33] {$3$};
            
            \node [rect node] (s41) [below=0cm of s31] {$3$};
            \node [rect node] (s42) [right=0cm of s41] {$2$};
            \node [rect node] (s43) [right=0cm of s42] {$1$};
            \node [rect node] (s44) [right=0cm of s43] {$4$};

            \node [empt node] (e1) [below=0cm of s44] {\phantom{M}};
            
        \end{tikzpicture}
    \end{subfigure}%
    \caption{Latin square of order $4$ generated by rule $150$.}
    \label{fig:r150-sq}
\end{figure}

\section{Mutually Orthogonal Cellular Automata (MOCA)}
\label{sec:mols-ca}

In this section we describe \emph{Mutually Orthogonal Cellular Automata (MOCA)}, which are families of bipermutive CA that form \emph{Mutually Orthogonal Latin Squares (MOLS)}. We start by illustrating the well-developed literature for linear MOCA, and then we move to the less explored nonlinear context. 

\subsection{The Linear Case and the Connection with Coprime Polynomials}
\label{subsec:lin-moca}

As shown in the previous section, any bipermutive CA gives rise to a Latin square when interpreted as an algebraic system. Thus, the next natural question is under which conditions the Latin squares generated by two bipermutive CA are orthogonal. Mariot et al.~\cite{mgfl-desi-2020} investigated this question, leveraging on the rich algebraic structure featured by linear CA. In what follows, we review the main characterization results obtained so far in this direction.

Recall from Section~\ref{subsec:ca} that the global rule of a linear CA can be described by a linear transformation whose transition matrix has the structure reported in Equation~\eqref{eq:ca-matr-f}. Let us consider now two bipermutive linear CA $F,G:\F_q^{2(d-1)} \to \F_q^{d-1}$ respectively defined by the coefficient vectors $(a_1,\cdots, a_d) \in \F_q^d$ and $(b_1,\ldots, b_d) \in \F_q^d$. Then, the Latin squares $C_F$ and $C_G$ are orthogonal if and only if the square matrix obtained by superposing the two transition matrices $M_F, M_G$ is invertible:

\begin{equation}
    \label{eq:syl-matr}
    M_{F,G} = 
    \begin{pmatrix}
        a_1    & \ldots & a_{d-1} & a_d & 0 & \ldots & \ldots & \ldots & \ldots & 0 \\
        0      & a_1    & \ldots  & a_{d-1} & a_d & 0 & \ldots & \ldots & \ldots & 0 \\
        \vdots & \vdots & \vdots & \ddots  & \vdots & \vdots & \vdots & \ddots & \vdots & \vdots \\
        0 & \ldots & \ldots & \ldots & \ldots & 0 & a_1 & \ldots & a_{d-1} & a_d \\
        b_1    & \ldots & b_{d-1} & b_d & 0 & \ldots & \ldots & \ldots & \ldots & 0 \\
        0      & b_1    & \ldots  & b_{d-1} & b_d & 0 & \ldots & \ldots & \ldots & 0 \\
        \vdots & \vdots & \vdots & \ddots  & \vdots & \vdots & \vdots & \ddots & \vdots & \vdots \\
        0 & \ldots & \ldots & \ldots & \ldots & 0 & b_1 & \ldots & b_{d-1} & b_d \\
    \end{pmatrix} \enspace .
\end{equation}
Indeed, from Section~\ref{subsec:cd} we know that $C_F$ and $C_G$ are orthogonal if and only if their superposition yields all ordered pairs in the Cartesian product $[N] \times [N]$ (where $N = q^{d-1}$) exactly once, or equivalently if the superposed mapping $\mathcal{H}: \F_q^{2(d-1)} \to \F_q^{(d-1)} \times \F_q^{(d-1)}$ defined for all $x \in \F_q^{2(d-1)} \times \F_q^{2(d-1)}$ as
\begin{equation}
\label{eq:sup-map}
    \mathcal{H}(x) = (F(x), G(x))
\end{equation}
is a bijective mapping. Since $M_{F,G}$ is precisely the matrix that defines $\mathcal{H}$, this mapping is bijective if and only if $M_{F,G}$ is invertible.

The authors of~\cite{mgfl-desi-2020} then observed that $M_{F,G}$ is the \emph{Sylvester matrix} of the two polynomials $p_f(X), p_g(X) \in \F_q[X]$ associated to the local rules of $F$ and $G$, and it is a well-known fact that the determinant of such matrix (also called the \emph{resultant}) is not null if and only if $p_f(X)$ and $p_g(X)$ do not have any factors in common~\cite{gelfand}. This observation resulted in the following characterization of orthogonal Latin squares defined by pairs of linear bipermutive CA:

\begin{theorem}[\cite{mgfl-desi-2020}]
\label{thm:lin-oca}
    Let $F,G:\F_q^{2(d-1)} \to \F_q^{d-1}$ be two linear bipermutive CA, with associated polynomials respectively $p_f,p_g \in \F_q[X]$. Then, the corresponding Latin squares $C_F$ and $C_G$ of order $N=q^{d-1}$ are orthogonal if and only if $\gcd(p_f,p_g) = 1$.
\end{theorem}
Theorem~\ref{thm:lin-oca} provides a simple and elegant condition to check whether the Latin squares induced by two linear bipermutive CA are orthogonal or not: it suffices to verify whether the polynomials associated to their local rules are relatively prime or not, which can be done efficiently with Euclid's algorithm.

Beside this characterization, the authors of~\cite{mgfl-desi-2020} also focused on two counting questions: namely, determining the number of orthogonal CA pairs and the size of the largest families of MOCA for a fixed diameter $d$. Concerning the first question, the problem is equivalent to counting the number of coprime polynomial pairs where both polynomials have degree $n=d-1$ and a nonzero constant term. While the counting question had already been settled for generic coprime pairs in several other papers~\cite{reifegerste00,benjamin07,hou09}, the variation where both polynomials must have a nonzero constant term required a separate treatment. In particular, using an approach based on recurrence equations, the authors of~\cite{mgfl-desi-2020} proved the following result:
\begin{theorem}[\cite{mgfl-desi-2020}]
\label{thm:count-lin-oca}
  Let $d>1$ and $n = d-1$. Then, the number of distinct pairs of coprime monic polynomials over $\F_q$ of degree $n$ and with a nonzero constant term (or equivalently, the number of orthogonal Latin square pairs generated by linear bipermutive CA of diameter $d$) is:
  \begin{equation}
    \label{eq:num-coprime}
    a_n = \frac{1}{2}\left(q (q-1)^3 \frac{q^{2n-2} - 1}{q^2 - 1} + (q-1)(q-2)\right) \enspace .
  \end{equation}
\end{theorem}
Interestingly, when $q=2$ (i.e. the binary CA case) Equation~\eqref{eq:num-coprime} becomes
\begin{equation}
    \label{eq:numb-cop-true}
    a_n = \frac{4^{n-1}-1}{3} \enspace ,
  \end{equation}
which is already known in the OEIS as the closed form of the recurrence equation for sequence A002450, that is connected to several other combinatorial and number-theoretic facts~\cite{a002450}.

Regarding the second question---determining the size of the largest families of MOCA for a given diameter $d$---the authors of~\cite{mgfl-desi-2020} considered an iterative construction to obtain a family of pairwise coprime polynomials of degree $n=d-1$ with a nonzero constant term. Intuitively, the idea underlying the construction is as follows: one starts by adding all irreducible polynomials of degree $n$ to the family, which clearly satisfy the desired condition of being pairwise coprime and having a nonzero constant term\footnote{Strictly speaking, the polynomial $X$ of degree 1 is also irreducible, but its constant term is null. However, this polynomial has been excluded from the construction in~\cite{mgfl-desi-2020}, since it does not correspond to a bipermutive CA (and thus cannot generate a Latin square).}. Then, for all $i \in \{1,\cdots \lfloor n/2 \rfloor\}$, one considers each irreducible polynomial of degree $i$ and multiplies it with a polynomial of degree $n-i$. The result is a polynomial of degree $n$ which, by construction, is coprime with all other polynomials already in the family, so it can be added. Finally, if $n$ is even all irreducible polynomials of degree $n/2$ are squared and added to the family. The following result shows that the size of the MOLS family resulting from this construction is optimal, i.e. no other linear bipermutive CA can be added that is orthogonal to all other members:
\begin{theorem}[\cite{mgfl-desi-2020}]
  \label{thm:char-max-mols}
  For any $n$ and $q$, the maximum number of MOCA of diameter
  $d = n+1$, or equivalently the maximum number of pairwise coprime monic
  polynomials of degree $n$ with nonzero constant term is:
  \begin{equation}
    \label{eq:char-max-mols}
    N_n = I_n + \sum_{k=1}^{\lfloor \frac{n}{2} \rfloor} I_k \enspace .
  \end{equation}
  where $I_k$ denotes the number of irreducible polynomials of degree $k$, computed through Gauss's formula~\cite{gauss-irr}.
\end{theorem}

A final question partially addressed in~\cite{mgfl-desi-2020} involves counting the number $T_n$ of MOCA families of maximal size $N_n$. In this regard, the authors proved a lower bound on $T_n$ by counting how many families can be obtained from their optimal construction (depending on the choices of the polynomials of degree $n-i$ at each step), and then proved that this lower bound is asymptotically close to the actual value of $T_n$.

Building upon the results of~\cite{mgfl-desi-2020}, Hammer and Lorch~\cite{Hammer23} considered the subclass of CA-based \emph{Factor Pair Latin Squares} (FPLS), a concept originally introduced in~\cite{Hammer17}. In general, a Latin square $L$ of order $N$ is a FPLS if, for each ordered factor pair $(a,b)$ such that $ab = N$, any tiling of size $a \times b$ in $L$ does not have repeated entries. The authors of~\cite{Hammer23} proved that linear CA defined over a prime alphabet generate factor pair Latin squares, and named them \emph{cellular factor pair Latin squares}. Further, the authors adapted Theorem~\ref{thm:char-max-mols} from~\cite{mgfl-desi-2020} to construct maximal families of mutually orthogonal cellular factor pair Latin squares, considering also the variant where the Latin squares are \emph{diagonal}, i.e., no symbols are repeated on the main diagonals of the squares. Interestingly, the authors showed that the size of these families is asymptotic to $p^n/n$ (with $p^n$ being the order of the Latin squares), improving on the bounds previously given for generic FPLS in~\cite{Hammer19}.

\subsection{The Nonlinear Case}
\label{subsec:nonlin-moca}

Cellular automata with nonlinear rules represent a significant departure from their linear counterparts, introducing greater complexity and unpredictability in their dynamics. This added complexity makes them particularly appealing as dynamical systems for cryptographic applications requiring high levels of security, such as pseudorandom number generators~\cite{formenti14,leporati13,leporati14}. On the other hand, the lack of linear relationships  makes the dynamics of nonlinear CA much more challenging to characterize. 

The interpretation of CA as algebraic systems considered in this paper is no exception: even if one is interested in investigating only the short-term behavior of a nonlinear CA, the inability to directly leverage the tools of linear algebra over finite fields is a significant obstacle. This situation is reflected in the scarcity of results obtained so far in the available literature concerning the construction of combinatorial designs via nonlinear CA.

From a theoretical standpoint, one of the main contributions in this field is represented by the work done in~\cite{mariot2017enumerating}, where the authors consider the enumeration of orthogonal Latin squares induced by nonlinear bipermutive CA over the binary alphabet $\F_2$. The core result presented in this paper is a necessary condition on the underlying local rules, which focuses on the property of \emph{pairwise balancedness}:

\begin{definition}[\cite{mariot2017enumerating}]
\label{def:pair-bal}
Two $d$-variable local rules $f,g: \F_2^d \rightarrow \F_2$ are
\emph{pairwise balanced} if the superposed function $(f,g): \F_2^d \rightarrow
\F_2^2$ defined as $(f,g)(x) = (f(x),g(x))$ is balanced, that is
$|(f,g)^{-1}(y_1,y_2)| = 2^{d-2}$ for all $(y_1,y_2) \in \F_2^2$.
\end{definition}

In other words, two local rules are pairwise balanced if the four pairs $(0, 0)$, $(1, 0)$, $(0, 1)$ and $(1, 1)$ occur an
equal number of times in the superposition of their truth tables. Then, the necessary condition is as follows:

\begin{lemma}[\cite{mariot2017enumerating}]
\label{lm:pwbal-bip}
Let $F,G: \F_2^{2(d-1)} \rightarrow \F_2^{d-1}$ be bipermutive CA respectively induced by local rules $f,g: \F_2^d \rightarrow \F_2$, and suppose that the associated Latin squares $C_F$ $C_G$ are orthogonal. Then, $f$ and $g$ are pairwise balanced.
\end{lemma}

Being a necessary condition, the above lemma gives at least a way to reduce the search space of candidate rule pairs whose CA give rise to orthogonal Latin squares. Moreover, the authors of~\cite{mariot2017enumerating} show a sufficient condition for two rules to be pairwise balanced. Recall from Equation~\eqref{eq:gen-func} in Section~\ref{subsec:ca} that a bipermutive rule of diameter $d$ can be defined in terms of its generating function computed on the central $d-2$ cells, since the leftmost and rightmost ones are simply XORed together. Then, the following result holds:

\begin{lemma}[\cite{mariot2017enumerating}]
\label{lm:gen-bip}
Let $f,g: \F_2^d \to \F_2$ be two bipermutive rules respectively defined by their $(d-2)$-variable generating functions $\varphi,\gamma: \F_2^{d-2} \to \F_2$. If $\varphi$ and $\gamma$ are pairwise balanced, then so are $f$ and $g$.
\end{lemma}

The authors then devise a combinatorial algorithm to exhaustively enumerate all orthogonal Latin squares generated by nonlinear bipermutive CA up to diameter $d=6$. Since Lemma~\ref{lm:gen-bip} is a sufficient condition, enumerating pairwise balanced generating functions would allow one to explore only a subset of the search space of interest. Therefore, the authors set out to give a combinatorial characterization to directly enumerate all pairwise balanced bipermutive rules. This characterization leverages an enumerative encoding for bipermutive rules originally introduced in~\cite{leporati13}, where the truth table is represented over a graph $G=(V,E)$. The set of vertices of this graph coincides with the space of all possible input vectors for the local rule, i.e. $V=\F_2^d$. Two vertices $v_1,v_2 \in V$ are connected by an edge if and only if $v_1$ and $v_2$ differ in either their leftmost or rightmost coordinate. This definition gives rise to a graph with $2^{d-2}$ connected components, where each component is composed of four vertices, and the degree of each vertex is 2.

The truth table of a bipermutive rule can be represented as a labeling $l_f: V \to \F_2$ over the vertices of such a graph: for each connected component, one assigns either a 0 or a 1 as a label to the node whose leftmost and rightmost coordinates are both 0. Due to the bipermutivity constraint, the labels of the the adjacent nodes must be set to the complementary value, while the label of the node where both the leftmost and rightmost coordinates are 1 must have the same label as the starting node. In this way, one can effectively enumerate the space of all bipermutive rules of diameter $d$ by enumerating all binary strings of length $d-2$.

To enumerate the space of pairwise balanced bipermutive rules, one can leverage the above graph-based representation by introducing a labeling function $l_{f,g}: V \to \F_2^2$ that assigns a \emph{pair} of bits to each node, instead of a single one. Beside the bipermutivity constraint for the single-rule labeling described earlier, which now applies component-wise, this new labeling must also satisfy the property that the four labels $(0,0)$, $(0,1)$, $(1,0)$ and $(1,1)$ occur an equal number $2^{d-2}$ of times. This turns out to be equivalent to the enumeration of balanced bistrings of length $2^{d-2}$, where a 0 and a 1 respectively identify a $(0,0)/(1,1)$ and a $(0,1)/(1,0)$ component. Further, since each type of connected component can be oriented in two different ways (i.e. assigning the first label either to the nodes that have the same or differing values on the leftmost and rightmost coordinates), this gives the following counting result:

\begin{lemma}[\cite{mariot2017enumerating}]
\label{lm:num-pwb}
The number of pairwise balanced pairs of bipermutive Boolean functions $f,g: \F_2^d \rightarrow \F_2$ of $d$ variables is:
\begin{equation}
\label{eq:count-pwval-bip}
\#Bal\mathcal{B}_d = \binom{2^{d-2}}{2^{d-3}} \cdot 2^{2^{d-2}} \enspace . 
\end{equation}
\end{lemma}

\begin{figure}[h]
\centering
\begin{tikzpicture}
[->,shorten >=1pt,auto,node distance=2cm,
  thick,main node/.style={circle,draw,font=\sffamily\bfseries},
        rep node/.style={circle,draw,font=\sffamily\bfseries}]

  \node[rep node] (1) {000};
  \node[above] at (1.north) {\bfseries{0,0}};
  \node[main node] (2) [below left of=1] {100};
  \node[above] at (2.north) {\bfseries{1,1}};
  \node[main node] (3) [below right of=2] {101};
  \node[above] at (3.north) {\bfseries{0,0}};
  \node[main node] (4) [below right of=1] {001};
  \node[above] at (4.north) {\bfseries{1,1}};
  
  \node[main node] (5) [right of=4] {110};
  \node[above] at (5.north) {\bfseries{1,0}};
  \node[main node] (6) [below right of=5] {111};
  \node[above] at (6.north) {\bfseries{0,1}};
  \node[main node] (7) [above right of=6] {011};
  \node[above] at (7.north) {\bfseries{1,0}};
  \node[rep node] (8) [above left of=7] {010};
  \node[above] at (8.north) {\bfseries{0,1}};
  
    \path[-]
    (1) edge node {} (2)
        edge node {} (4)
    (2) edge node {} (3)
        edge node {} (1)
    (3) edge node {} (4)
        edge node {} (2)
    (4) edge node {} (1)
        edge node {} (3);
    \path[-]
    (5) edge node {} (6)
        edge node {} (8)
    (6) edge node {} (7)
        edge node {} (5)
    (7) edge node {} (8)
        edge node {} (6)
    (8) edge node {} (5)
        edge node {} (7);
\end{tikzpicture}
\caption[]{Graph representation of the pairwise balanced pair of diameter $d=3$ composed of rules 90 and 150.}
\label{fig:pwbal-graph}
\end{figure}
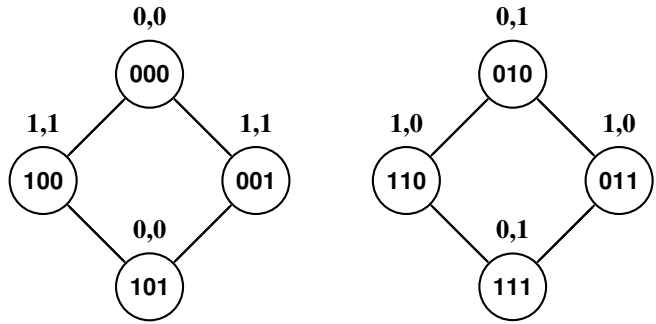

As an example, Figure~\ref{fig:pwbal-graph} depicts the graph representation of the pairwise balanced pair composed of rules 90 and 150 of diameter $d=3$ (which, incidentally, also gives rise to a pair of orthogonal CA).

The combinatorial algorithm conceived by the authors of~\cite{mariot2017enumerating} to enumerate all pairwise balanced pairs up to $d=6$ is then a variation of Knuth's algorithm to generate all bitstrings of a fixed Hamming weight~\cite{knuth2022}. The results report the number of orthogonal Latin squares generated by the corresponding bipermutive CA, along with the distribution of nonlinearity values of the underlying local rules.

Clearly, exhaustive enumeration becomes unfeasible beyond diameter $d=6$: following Equation~\eqref{eq:count-pwval-bip}, the number of pairwise balanced bipermutive rules is already approximately $2.58 \cdot 10^{18}$ for $d=7$. To cope with this problem, the authors of \cite{mariot2017evolutionary} considered the construction of orthogonal Latin squares generated by bipermutive CA as an optimization problem, to be tackled with Genetic Algorithms (GA) and Genetic Programming (GP).

To address the challenge of optimizing the evolution of CA-based Orthogonal Latin Squares (OLS), the authors exploited the property that bipermutive rules involving $d$ variables can be described using generating functions of $d-2$ variables. Consequently, the genotype of each individual in the population encodes a pair of generating functions, ensuring that the resulting phenotype corresponds to a pair of bipermutive CA capable of generating two Latin squares. As a result, the optimization process in both GA and GP can concentrate on enhancing the orthogonality and nonlinearity of the solutions, bypassing the need to verify the row-column permutation property.

For GA, the authors proposed three distinct encodings for the chromosomes representing candidate solutions. The first encoding combines the truth tables of the two generating functions into a single bitstring, allowing standard crossover and mutation operators to be directly applied. In the second encoding, the two generating functions are treated independently, with genetic operators applied separately to each function. Lastly, the approach takes advantage of two empirical observations derived from the search experiments of~\cite{mariot2017enumerating}, namely pairwise balancedness of the local rules, and the sufficient condition of Lemma~\ref{lm:gen-bip} that pairwise balanced generating functions induce pairwise balanced rules. Leveraging these insights, the authors designed the third encoding, where the genotype is structured as a balanced quaternary string. To ensure the balancedness property is maintained in quaternary strings, the authors develop specialized crossover and mutation operators tailored for GA~\cite{manzoni20}.

In contrast, for GP, the authors employed an encoding similar to the GA's double bitstring representation, where each generating function is modeled as an independent Boolean tree. Both GA and GP are evaluated on the smallest problem instances beyond the reach of exhaustive search, specifically the sets of generating function pairs with $5$ and $6$ variables. These instances correspond, respectively, to the sets of bipermutive CA pairs with local rules of $n = 7$ and $n = 8$ variables, or equivalently to the sets of CA-based Latin squares pairs of size $64\times 64$ and $128\times 128$. In the experiments, the authors adopted two fitness functions to be minimized: the first only counts the number of repeated pairs in the superposition of two Latin squares, while the second also adds a penalty if either or both generating functions are linear.

The main insight obtained from the experiments is that GP outperforms by far GA in evolving orthogonal CA. However, if the fitness function only focuses on the orthogonality property, GP interestingly converges almost always to linear pairs, probably because they have a simpler description in terms of Boolean trees. Forcing the additional constraint of nonlinearity in the fitness function ensures that GP always finds pairs of nonlinear bipermutive CA that generate orthogonal Latin squares.

\section{Applications to Cryptography}
\label{sec:appl}

In this section, we survey some applications to cryptography of the results reviewed in the previous section. In particular, we consider how MOCA can be used to define threshold secret sharing schemes, pseudorandom number generators, and Boolean functions with significant cryptographic properties.
	
\subsection{Threshold Secret Sharing Schemes}
\label{subsec:sss-ac}
Recall from Section~\ref{subsec:cd} that a $(t,n)$ perfect threshold secret sharing scheme enables a dealer to share a secret $S$ with a set of $n$ participants, such that any subset of at least $t$ players can uniquely reconstruct $S$ by combining together the respective shares distributed by the dealer. On the other hand, any subset of less than $t$ players gain no information on the value of $S$, meaning that by knowing at most $t-1$ shares, any value of the secret is equally likely.

A sets of $n$ MOLS $L_1,\ldots, L_n$ of order $N$ is equivalent to a $(2, n)$ threshold scheme, where both the secret and the shares range in the set $[N] = \{1,\cdots, N\}$. The idea of the protocol is the following: given the secret $S \in [N]$, the dealer uses it to index a \emph{row} of the $n$ MOLS. Then, the dealer samples with uniform probability a \emph{column} $R \in [N]$. The entries of the $n$ squares at the coordinates $(S, R)$ correspond to the shares, i.e. for all $i \in \{1,\cdots, n\}$ the $i$-th player receive the value $B_i = L_i(S, R)$ as its share. Further, the dealer also publishes the set of $n$ MOLS, so it is known to the players. Later, if two players $i,j$ want to retrieve the value of $S$, they pool together their shares, thus obtaining the pair $(B_i, B_j) = (L_i(S, R), L_j(S, R))$.

Since all Latin squares in the MOLS family are pairwise orthogonal, $i$ and $j$ can check where the unique entry $(L_i(S, R), L_j(S, R)$ occurs in the superposition of $L_i$ and $L_j$. The row coordinate of this entry corresponds to the original secret $S$. On the other hand, if---without loss of generality---only the $i$-th player attempts to reconstruct the secret, its share $B_i$ occurs exactly once in each row and in each column of the Latin square $L_i$. Therefore, if the column $R$ has been chosen uniformly at random by the dealer, any row (and thus any value of the secret) is equally likely. Hence, any single player does not gain any information about the value of $S$ by just knowing its own share.

Given the results reviewed in Section~\ref{sec:mols-ca}, any set of $n$ MOCA can be used to implement a $(2, n)$ threshold scheme following the above protocol. However, in practical settings the size of the secret is usually too large to grant the explicit publication of the squares in a MOLS family. For instance, if a dealer wants to share a 128-bit key with the players, then the order of the corresponding Latin squares would be $2^{128}$. Consequently, it is necessary to find an \emph{implicit} efficient representation of the Latin squares, both to sample their entries and to check in retrospect the row and column coordinates of a superposed pair of entries. In the case of linear MOCA, this compact representation can be obtained by encoding each Latin square just by the coefficients of the underlying linear CA local rule. The computation of the shares by the dealer amounts to initializing the left half of the CA to the value of the secret $S$, while the right half is the random value $R$ used to index the column of the squares. Finally, the output vector obtained by evaluating the NBCA over this initial configuration is the share to be given to the player. On the other hand, recovering the secret from a pair of output vectors $(y_i, y_j)$ held by players $i,j$ can be obtained by constructing the Sylvester matrix associated to the polynomials of the CA $F_i,F_j$. Then, the players compute the inverse of this matrix and multiply it with the concatenation of the two vectors $y_i||y_j$. The result is the initial configuration of the CA, whose left half corresponds to the secret $S$.

More in general, Mariot and Leporati explored the inversion problem for MOCA in \cite{ml-acri-2018}. There, the authors presented an algorithm to invert a pair of output configurations obtained by evaluating two orthogonal CA $F,G$ on the same input vector. The algorithm is based on the \emph{coupled de Bruijn graph}, which is obtained by superposing the edge labels of the de Bruijn graphs associated to the two CA. Given a pair of output configurations $(w_1, \ldots, w_{d-1})$ and $(z_1, \ldots, z_{d-1})$ respectively computed by $F$ and $G$, the algorithm starts by looking up all edges on the coupled de Bruijn graph labeled as $(w_1, z_1)$. Then, for each of these nodes it applies a Depth-First Search (DFS) strategy to follow the edge paths subsequently labeled as $(w_2, z_2)$, $(w_3, z_3)$, and so on. Due to the orthogonality of the two CA, only one of these DFS visits succeeds in finding a complete edge path labeled from $(w_1, z_1)$ up to $(z_{d-1}, z_{d-1})$. By applying the fusion operator on the sequence of vertices visited by this path, one obtains the initial configuration of the CA.

The advantage of the above procedure is that it can be applied to the inversion of any pair of orthogonal CA, independently from the linearity of their local rules. However, an obvious drawback is that the size of the coupled de Bruijn graph is still exponential in the diameter of the CA, so in general it does not yield an efficient representation. The authors argue that this shortcoming can be addressed by modifying the inversion algorithm to leverage the \emph{Algebraic Normal Form} (ANF)~\cite{carlet21} of the local rules and their permutivity properties. In particular, one can consider pairs of orthogonal CA whose local rules have a low algebraic degree to obtain a compact representation for the inversion algorithm.

\subsection{Pseudorandom Number Generators}
\label{subsec:prng}

CA have long been considered as Pseudorandom Number Generators (PRNGs) for cryptographic purposes. As mentioned in the Introduction, the very first work by Wolfram~\cite{wolfram85} that applied CA to cryptography was precisely about a PRNG based on the dynamics of a CA, which produced the keystream sequences of a Vernam-like stream cipher. The idea was to initialize the seed of the PRNG as the initial configuration of a one-dimensional periodic CA equipped with the elementary rule 30, and then iterate the CA for as many steps as the length of the desired keystream. The trace of the central cell was thus taken as the keystream to be XORed with the plaintext. Albeit simple, the idea turned out to be vulnerable to a few critical attacks~\cite{meier-staff,koc-ca}. Beside this, Daemen et al.~\cite{daemen94-1} observed that such CA-based generators inevitably suffer from poor diffusion, since the speed a small difference in the input configuration can propagate with is always bound by the diameter of the local rule. Hence, one needs to iterate the CA for a considerable number of steps to ensure that the difference has spread in the whole lattice.

To address this problem, Mariot~\cite{m-naco-2023} investigated the use of orthogonal CA to generate pseudorandom sequences. The main motivation is that a pair of orthogonal Latin squares defines a $(2, 2)$-\emph{multipermutation}~\cite{vaudenay94}. Given two orthogonal CA $F,G: \F_q^{2(d-1)} \to \F_q^{d-1}$, this means that for any two distinct pairs of $(d-1)$-cell configurations $(x,y), (x',y') \in \F_q^{d-1} \times \F_q^{d-1}$ the corresponding input/output 4-tuples $(x, y, F(x, y), G(x, y))$ and $(x', y', F(x', y'), G(x', y'))$ cannot agree on more than one coordinate. This ensures that the superposed permutation $\mathcal{H}$ defined in Equation~\eqref{eq:sup-map} has an optimal diffusion in terms of $(d-1)$-cell blocks.

However, orthogonal CA are meant to be algebraic systems, and thus they can be iterated only for a limited number of steps (usually, a single one). On the other hand, a PRNG has to stretch an initial short seed into an arbitrarily long pseudorandom sequence. Hence, the author of~\cite{m-naco-2023} defines a dynamical system by iterating the superposed permutation $\mathcal{H}$. In other words, the PRNG is initialized with a seed of length $2(d-1)$. The two orthogonal CA are applied to this seed, after which one obtains a pair of configurations of length $(d-1)$. This two configurations are then concatenated to obtain a new vector of length $2(d-1)$, to which the orthogonal CA can be again applied, and so on.

Since the map $\mathcal{H}$ is a permutation, the dynamics of this system is reversible, or equivalently it can be decomposed into cycles. The work in~\cite{m-naco-2023} proceeds by studying the cycle decomposition of this permutation, focusing on both nonlinear and linear orthogonal CA. In the former case, the author provides some empirical results on the distribution of maximum cycle lengths for orthogonal CA of diameters between 4 and 8, remarking that the largest cycles of size $2^{2(d-1)}-1$ occurs only for linear pairs. Then, the focus narrows on linear orthogonal CA, where the author remarks that the maximum cycle length of the dynamical system corresponds to the order of the Sylvester matrix associated to the CA. This allows to leverage the theory of Linear Modular Systems (LMS) to give a complete characterization of the cycle structure of linear orthogonal CA. In particular~\cite{m-naco-2023} proves that the its maximal period of $q^{2n} - 1$ is attained when the Sylvester matrix $M_{f,g}$ has a primitive minimal polynomial. Extensive enumeration experiments up to $d = 16$ (for $q=2$) and $d = 13$ (for $q=3$) are then carried out to generate all linear OCA pairs achieving the maximum cycle length.

\subsection{Boolean Functions}
\label{subsec:boolfun}
Boolean functions play a crucial role in symmetric cryptography, for instance, in the design of stream ciphers based on the combiner and the filter model~\cite{carlet21}. In the combiner model, a Boolean function $f: \F_2^n \to \F_2$ is used to compute the next keystream bit by combining the outputs of $n$ Linear Feedback Shift Registers (LFSRs), while in the filter model $f$ taps the registers of a single LFSR of size $n$. In both cases, the security of the stream cipher can be reduced to the \emph{cryptographic properties} of the underlying Boolean function. The rationale is that if a Boolean function does not satisfy a particular property, the cipher is vulnerable to a specific attack.

\subsubsection{Bent functions}
One of the most important cryptographic criteria for a Boolean function $f: \F_2^n \to \F_2$ is \emph{nonlinearity}, which can be defined as the Hamming distance of the truth table of $f$ from the set of all linear functions $L_a(x) = a \cdot x = a_1x_1 \oplus \ldots \oplus a_n x_n$. In stream ciphers, functions with low nonlinearity are vulnerable to fast correlation attacks~\cite{meier88}. Thus, as a design criterion, a Boolean function $f: \F_2^n \to \F_2$ used in the combiner or filter model should lie at a high Hamming distance from all $n$-variable linear functions.

The most convenient way to compute the nonlinearity of a Boolean function $f: \F_2^n \to \F_2$ is the \emph{Walsh transform} $W_f: \F_2^n \to \Z$, which is defined as:
\begin{equation}
    \label{eq:wht}
    W_f(a) = \sum_{x \in \F_2^n} (-1)^{f(x) \oplus a \cdot x} \enspace ,
\end{equation}
for all $a \in \F_2^n$. In essence, the coefficient $W_f(a)$ measures the correlation of $f$ with the linear function $a\cdot x$. Then, the nonlinearity of $f$ equals:
\begin{equation}
    \label{eq:nl}
    nl(f) = 2^{n-1} - \frac{1}{2}\max_{a \in \F_2^n} \{|W_f(a)|\} \enspace .
\end{equation}
Ideally, a Boolean function should have minimum correlation with every linear function to achieve a high nonlinearity. However, \emph{Parseval identity} states that the sum of the squared Walsh coefficients equals $2^{2n}$ for every $n$-variable Boolean function. Consequently, the nonlinearity of a Boolean function $f: \F_2^n \to \F_2$ is bounded above by the so-called \emph{covering radius bound}:
\begin{equation}
    \label{eq:crb}
    nl(f) \le 2^{n-1} - 2^{\frac{n}{2}-1} \enspace .
\end{equation}
Functions that satisfy the above bound with equality are called \emph{bent}, and they exist only for even values of $n$, since the Walsh transform yields only integer values.

Bent functions are connected to several other combinatorial designs, including Hadamard matrices~\cite{rothaus76}, strongly regular graphs~\cite{bernasconi01} and partial spreads~\cite{dillon74}. Several algebraic constructions of bent functions leverage some of these connections~\cite{carlet21,mesnager16}. Further, there exists a construction of Hadamard matrices from families of MOLS, due to Bush~\cite{bush73}.

Starting from this motivation, Gadouleau et al.~\cite{gmp-iacr-2020} proved that the construction of mutually orthogonal Latin squares from linear bipermutive CA in~\cite{mgfl-desi-2020} gives rise to a Hadamard matrix which is in turn associated to a bent function. The authors then addressed the related existence and counting questions, proving that these functions exist only if the degree of the underlying coprime polynomials is either 1 or 2. Incidentally, the counting results settled also the question of determining the exact value of $T_n$ discussed in Section~\ref{subsec:lin-moca} for $n=2$. Further, the authors showed that the functions resulting from this construction belong to the partial spread class $\mathcal{PS}$. The construction has been later simplified in~\cite{gmp-desi-2023} by showing a direct way to construct partial spreads from linear recurring sequences, without passing through the characterization steps with Latin squares, linear CA and Hadamard matrices. Additionally, the authors performed a computational analysis of the ranks of all bent functions of $n=6$ and $n=8$ variables generated by this construction. The main finding is that bent functions obtained from MOCA families whose underlying polynomials have degree 1 coincide with the functions in the $\mathcal{PS}_{ap}$ class (also called the \emph{Desarguesian spread}). On the other hand, for degree 2 the analysis of the ranks found functions that are equivalent neither to the Maiorana-McFarland nor to the Desarguesian spread classes.

\subsubsection{Correlation Immunity}
Another important cryptographic property used in stream ciphers is \emph{correlation immunity}. Formally, a Boolean function $f: \F_2^n \to \F_2$ is $t$-th order correlation immune if the probability distribution of its output remains unchanged when fixing at most $t$ input variables. In the context of the combiner model, functions that are not correlation immune of a sufficiently high order are vulnerable to correlation attacks~\cite{siegenthaler85}. Moreover, correlation immune functions are also useful in the design of \emph{masking countermeasures} against side-channel attacks~\cite{carlet12}.

Camion et al.~\cite{camion91} proved a characterization of correlation immune functions in terms of orthogonal arrays. Namely, a Boolean function is $t$-th order correlation immune if and only if its support $supp(f) = \{x \in \F_2^n: f(x) \neq 0\}$ is a binary orthogonal array of strength $t$. The paper~\cite{mariot2023building} by Mariot and Manzoni addresses the challenge of designing correlation-immune Boolean functions through MOCA, leveraging the above orthogonal array characterization. As mentioned in Section~\ref{subsec:cd}, a family of $k$ MOLS of order $s$ is equivalent to an orthogonal array of strength $2$ with entries in $[s] = \{1,\ldots, s\}$. Thus, the authors of~\cite{mariot2023building} considered the OA arising from MOCA families, encoding each entry as a binary vector. The main theoretical result is that the resulting array is indeed a binary OA of strength at least $2$. Therefore, families of MOCA define the support of correlation immune functions of order at least 2. Finally, the authors of~\cite{mariot2023building} performed an exhaustive search of all families of 3 MOCA defined by local rules of diameter 4 and 5, observing that all resulting functions are actually correlation immune of order at least 3.

\section{Conclusions and Future Directions}
\label{sec:outro}

In this survey, we explored the interplay between combinatorial designs and cellular automata, with a particular focus on their algebraic interpretation and cryptographic applications. We reviewed key results concerning the construction of mutually orthogonal Latin squares via bipermutive CA and examined the implications of these structures for secret sharing schemes, pseudorandom number generators, and Boolean functions.

Our discussion began by presenting the algebraic perspective on CA, where the short-term behavior of certain bipermutive CA can be associated with quasigroups and Latin squares. We then examined how families of mutually orthogonal CA can be systematically constructed, particularly in the case of linear CA, where the orthogonality condition corresponds to the coprimality of their associated polynomials. Although significant progress has been made in the linear case, the nonlinear setting remains an open research direction, with recent approaches leveraging combinatorial algorithms and optimization techniques to construct nonlinear MOCA.

Beyond theory, we surveyed applications of MOCA to the design of cryptographic primitives. We highlighted how CA-based MOLS can be efficiently used to implement $(2, n)$ threshold secret sharing schemes, offering a compact representation of the involved Latin squares. Additionally, we discussed the use of CA in generating bent and correlation immune Boolean functions with desirable cryptographic properties, and their potential in the generation of pseudorandom sequences.

Despite these advances, several open problems remain. We conclude the paper by highlighting a few directions that we deem particularly interesting and worth exploring in future research.

\subsection{Improve Lower Bounds on the Number of MOCA Families}
One key challenge is obtaining better lower bounds on $T_n$, the number of maximal MOCA families of a given size, or achieving a complete combinatorial characterization of this quantity. While a lower bound for $T_n$ has been established through a constructive approach, determining the exact number of distinct maximal MOCA families remains an open question. A deeper understanding of the combinatorial structure of MOCA could lead to more efficient construction techniques and potentially new cryptographic applications. This problem shares similarities with the enumeration of bent functions, particularly in the context of combinatorial constructions. The work in~\cite{gmp-desi-2023} provided the exact values for $T_1$ and $T_2$. Further developing these techniques might give new insights into the enumeration of maximal MOCA families when polynomials of a larger degrees are involved.

\subsection{Generalization to Hypercubes and Orthogonal Arrays}
Another promising direction is the generalization of CA-based Latin squares to Latin hypercubes. A Latin hypercube of dimension $k$ and order $n$ is a $k$-dimensional array where each row, column, and higher-dimensional analog contains each symbol exactly once. A preliminary investigation of Latin hypercubes defined by CA has already been carried out in~\cite{gm-automata-2020}, where the authors proposed a construction based on linear CA. However, there are still interesting open problems in this research thread, such as:
\begin{itemize}
    \item \emph{Finding sufficient conditions for orthogonality}. Unlike the case of Latin squares, where the orthogonality condition is well understood for CA-based constructions (i.e., via coprimality of polynomials in the linear case), the corresponding condition for mutually orthogonal subsets of Latin hypercubes is an open problem. A precise characterization of this property derived from CA would enable the systematic generation of such structures. A starting point could be the generalized resultant for families of $k$ polynomials introduced by Dei{\ss}ler in~\cite{deissler}, since it has the most similar structure to the superposition of the transition matrices of $k$ linear CA.
    \item \emph{Secret sharing schemes with threshold larger than 2}. The connection between MOLS and $(2, n)$-threshold secret sharing schemes is well established, but Latin hypercubes could provide a natural extension to threshold schemes where any $t$ out of $n$ shares can reconstruct the secret. This would generalize the current CA-based $(2, n)$ schemes to more general $(t, n)$ threshold sharing schemes. Moreover, since $(t, n)$ schemes are a particular form of MDS codes, one could consider the construction MDS matrices with CA as a further application, which would be relevant in the design of diffusion layers for block ciphers~\cite{gupta2019cryptographically}.
\end{itemize}
    
\subsection{Improve Results on Nonlinear MOCA}
Another significant open question is the characterization of nonlinear MOCA in terms of systems of nonlinear equations, which would allow the application of computational algebra techniques to solve them. Unlike the linear case, where the orthogonality condition corresponds to polynomial coprimality, nonlinear CA introduce additional complexity, as their evolution follows nonlinear transformations that are difficult to analyze directly. To formalize this problem, consider a set of $k$ nonlinear bipermutive CA $F_1,F_2,\dots,F_k$ over the finite field $\F_q$. The requirement that the associated Latin squares are mutually orthogonal translates into a set of nonlinear constraints on the local rules $f_1,f_2,\dots,f_k$ that can be expressed as a set of algebraic equations over $\F_q.$ Specifically, the function $H(x)=(F_1(x),F_2(x),\dots,F_k(x))$ must be a bijection, meaning that each output tuple appears exactly once across all possible inputs. This condition leads to a system of nonlinear polynomial equations in multiple variables. A natural approach for solving these equations is through Gr\"obner bases, a fundamental tool from computational algebra that enables the simplification and solution of polynomial systems: in particular, Gr\"obner bases can be used to determine whether the system has solutions, to count the number of solutions, or to construct explicit solutions.

The application of Gr\"obner bases to the MOCA problem could also be compared to similar approaches in the study of nonlinear Boolean functions with optimal cryptographic properties. For example, in Boolean function design, Gr\"obner bases have been used to analyze and construct functions with maximum algebraic immunity~\cite{cid2009block}, a key property for resisting algebraic attacks on stream ciphers. Ultimately, characterizing nonlinear MOCA through systems of polynomial equations would bridge the gap between combinatorial design theory and computational algebra, opening new pathways for constructing and analyzing CA-based Latin squares.

\subsection*{Acknowledgements}
Giuliamaria Menara was supported by the PRIN 2022 PNRR project entitled ``Cellular Automata Synthesis for Cryptography Applications (CASCA)" (P2022MPFRT) funded by the European Union – Next Generation EU. Luca Manzoni was supported by the Horizon 2020 MSCA project entitled ``Application-driven Challenges for Automata Networks and Complex Systems (ACANCOS)'' (MSCA-SE-101131549) funded by the European union under the HORIZON-TMA-MSCA-SE action.

\bibliographystyle{abbrv}
\bibliography{references}

\end{document}